# On embedded minimal surfaces of Costa-Hoffman-Meeks type in hyperbolic space.

25th May 2017

Asun Jiménez Grande[1], Graham Smith[2]

**Abstract:** We present a new construction of embedded minimal surfaces in hyperbolic space with 3 asymptotically totally geodesic ends and arbitrary finite genus.

**Key Words:** Minimal surfaces, hyperbolic space, desingularisation, Costa-Hoffman-Meeks surfaces

**AMS Subject Classification:** 53A10

---

[1] Instituto de Matemática, UFF, Rua Professor Marcos Waldemar de Freitas Reis, Bloco H - Campus do Gragoatá, São Domingos, 24.210-201, Niterói, RJ - BRAZIL

[2] Instituto de Matemática, UFRJ, Av. Athos da Silveira Ramos 149, Centro de Tecnologia - Bloco C, Cidade Universitária - Ilha do Fundão, Caixa Postal 68530, 21941-909, Rio de Janeiro, RJ - BRAZIL

# 1 - Introduction.

**1.1 - Introduction and main result.** Ever since the pioneering work [8], [9] & [10] of Kapouleas, desingularisation constructions have provided a rich source of minimal and constant mean curvature surfaces with a variety of geometric properties. Typically, geometrically complex surfaces are obtained by splicing elementary components, such as planes and catenoids, along singly periodic Scherk surfaces of high genus. In certain cases, however, the splice may also be performed along other types of surface, which usually provides a finer control over the geometries of the surfaces constructed. In this paper, we will be interested in the case where the splice is performed along a Costa-Hoffman-Meeks surface (c.f. [3] & [7]).

Costa-Hoffman-Meeks surfaces have already been used in various desingularisation constructions, exemplified, perhaps, by the work [6] of Hauswirth & Pacard and [11] of Mazzeo & Pacard. However, all such constructions have to date relied upon an ingenious argument based on the geometry near infinity of a certain intermediate surface which, significantly, breaks down when this geometry is not of the correct type. Indeed, denoting the intermediate surface by $\Sigma$ and its Jacobi operator by $J$, the main analytic challenge in any desingularisation argument lies in constructing a suitable right inverse for $J$ over $\Sigma$. In the above-cited papers, it turns out that $\Sigma$ and $J$ are respectively - at least intrinsically - close to a Costa-Hoffman-Meeks surface $\Sigma_g$ and its Jacobi operator $J_g$, and the right inverse of $J$ is then readily obtained by a perturbation argument (c.f. [13] & [14]). In cases where such an approximation is no longer valid, a different and more arduous approach is required. One such case is that of properly embedded translating solitons of the mean curvature flow in $\mathbb{R}^3$, which are minimal surfaces with respect to a certain conformal metric, and which were studied in detail by the second author in [18]. Another is that of properly embedded minimal surfaces of finite topological type in three-dimensional hyperbolic space $\mathbb{H}^3$, and it is these that will form the subject of this paper.

A number of constructions of properly embedded minimal surfaces of finite topological type in $\mathbb{H}^3$ are already known. However, since hyperbolic space is essentially much "bigger" than $\mathbb{R}^3$, a classification of such surfaces as straightforward as that established for the Euclidean case by Ossermann in [16] is not to be expected. Indeed, already in the stable case, in contrast to the uniqueness results [2] of do Carmo & Peng and [4] of Fischer-Colbrie & Schoen, Anderson shows in [1] that there exists a properly embedded, simply connected, stable minimal surface bounded by every Jordan curve in the ideal boundary of $\mathbb{H}^3$. More generally, using variational techniques, Oliveira & Soret construct in [15] properly embedded stable minimal surfaces in $\mathbb{H}^3$ of arbitrary finite topology. Finally, by means of a bridging argument, Martín & White construct in [12] another family of properly embedded minimal surfaces in $\mathbb{H}^3$, also of arbitrary finite topology, but not clearly related to that of Oliveria & Soret.

In this paper, we complement the works of Oliveira & Soret and Martín & White by constructing what we consider to be the closest hyperbolic analogues of the Costa-Hoffman-Meeks surfaces. In particular, these surfaces, whose geometric structures we are able to describe with a fair degree of precision, are quite distinct from the ones constructed above.

In order to state the theorem, we first recall some basic geometric properties of the Costa-Hoffman-Meeks surfaces (c.f. [7]). For every positive integer $g$, the Costa-Hoffman-Meeks surface $\Sigma_g$ is a properly embedded minimal surface in $\mathbb{R}^3$ of genus $g$ with 3 ends, each of which we assume to be a graph over an unbounded annulus in $\mathbb{R}^2$. For all $0 \leq k \leq g$, this surface is invariant under reflection in the plane containing the $z$-axis which makes an angle of $k\pi/(g+1)$ with the $x$-axis. We refer to the group of symmetries of $\mathbb{R}^3$ generated by these reflections as the group of *horizontal symmetries* of $\Sigma_g$.* We prove

**Theorem A**

*Fix $g \in \mathbb{N}$ and $\eta \ll 1$. For all sufficiently large $\Lambda > 0$, and for all $\epsilon > 0$ and $R > 0$ satisfying*

$$\epsilon R^{5-2\eta} \leq \frac{1}{\Lambda} \text{ and } \epsilon R^{5-\eta} \geq \Lambda, \tag{1.1}$$

*there exists a complete, embedded minimal surface $\Sigma$ in $\mathbb{H}^3$ of genus $g$ with 3 ends. Furthermore*

*(1) $\Sigma$ is preserved by the group of horizontal symmetries of the Costa-Hoffman-Meeks surface $\Sigma_g$;*

*(2) $\Sigma \setminus B(\epsilon R)$ consists of three disjoint ends each of which converges towards the same horizontal, totally geodesic plane as $\Lambda$ tends to infinity; and*

*(3) Upon rescaling by a factor of $1/\epsilon$, $\Sigma \cap B(2\epsilon R)$ converges towards the Costa-Hoffman-Meeks surface $\Sigma_g$ as $\Lambda$ tends to infinity.*

**Remark:** All notation and terminology used in this paper is explained in detail in Appendix A.

**Remark:** Theorem **A** follows from Theorems 5.1.4 and 5.1.5, below.

**Remark:** The quantity $\Lambda$ serves to drive $\epsilon$ to zero and $R$ to infinity in a controlled manner: heuristically, $\epsilon$ is a little bit smaller than $R^{-5}$. The quantity $\epsilon$ determines the scaling factor of the Costa-Hoffman-Meeks surface. In particular, for sufficiently large $\Lambda$, distinct values of $\epsilon$ should yield distinct minimal surfaces in $\mathbb{H}^3$. Finally, the quantity $R$ determines how far along the end the glueing operation is carried out. Once $\epsilon$ has been fixed, we see no reason to expect that different values of $R$ should yield different minimal surfaces.

**1.2 - Techniques.** The proof of Theorem A essentially follows the standard desingularization argument first laid out by Kapouleas in [8], [9] & [10]. It is summarised as follows. First, using Fermi coordinates about a complete geodesic, $\mathbb{H}^3$ is identified with $\mathbb{R}^3$ furnished with the metric

$$g = dr^2 + \sinh^2(r)d\theta^2 + \cosh^2(r)dt^2.$$

With $\Lambda$, $\epsilon$ and $R$ as above, the Costa-Hoffman-Meeks surface $\Sigma_g$ is then rescaled by a factor of $\epsilon$ about the origin. In particular, the intersection of this rescaled surface with the annular prism $A(\epsilon R, \infty) \times \mathbb{R}$ consists of three disjoint ends, each of which is a graph over $A(\epsilon R, \infty)$. Using cut-off functions, each of these ends is then modified so as to coincide outside the solid cylinder $B(2\epsilon R) \times \mathbb{R}$ with a rotationally symmetric minimal end in $\mathbb{H}^3$ with parameters carefully chosen so as to minimise the curvature of the modified surface. A fixed-point argument then allows us to conclude by perturbing this surface into one which is minimal in $\mathbb{H}^3$.

---

* The complete symmetry group of $\Sigma_g$ is the dihedral group generated by the elements $A$ and $B$, where $A$ is reflection in the $(x-z)$-plane and $B$ is rotation by an angle of $k\pi/(g+1)$ about the $z$-axis followed by reflection in the $(x-y)$-plane (c.f. [7]).

It is the final stage of this process, which requires the construction of the right inverse of the Jacobi operator of the modified surface, that involves by far the most work. Here we closely follow the analysis of the Costa-Hoffman-Meeks surface already carried out by the second author in [18]. In particular, we use the detailed estimates already derived in Section 6.2 of that paper.

Finally, we believe that two main features of this analysis are worthy of mention. Firstly, in order for the estimates to decay correctly, we require that $\epsilon$ tend to zero at some rate between $R^{-4}$ and $R^{-5}$, and it is this condition, established by lengthy experimentarion with a variety of alternatives, which motivates the formulation of (1.1). Secondly, in contrast to desingularisation constructions which use singly periodic Scherk surfaces, the symmetries of the problem barely contribute to the derivation of the estimates that we use. Instead, we obtain these estimates, in part through a more careful and thorough analysis, but also through the introduction, in Section 2.3, of what we chose to call the *hybrid norm*

$$\|f\|_{m,\alpha} = \|f\|_{C^{k,\alpha}(\mathbb{H}^2)} + \frac{1}{(\epsilon R)} \|f\|_{H^k(\mathbb{H}^2)}.$$

It turns out that this norm, which is a weighted combination of the usual Hölder and Sobolev norms, with weight depending on $\epsilon$ and $R$, succinctly encapsulates the singular behaviour of this construction as $\Lambda$ becomes large, thus permitting us to complete the construction, even when the symmetry group is small.

**1.3 - Overview.** The paper is organized as follows.

In Section 2, we study the analytic properties of rotationally symmetric minimal ends in hyperbolic space, proving, in particular, the invertibility of their Jacobi operators with respect to suitably weighted hybrid norms.

In Section 3, we present the geometric construction, describing in detail the smooth families of candidate surfaces out of which the minimal surface will be selected.

In Section 4, using the estimates derived in Section 6.2 of [18], we construct right inverses of the Jacobi operators of the surfaces constructed in Section 3, together with estimates of their norms.

Finally, in Section 5, the Schauder fixed-point theorem is applied to prove the existence of an embedded minimal surface within the family constructed in Section 3.

## 2 - Catenoidal ends in hyperbolic space.

**2.1 - Surfaces of revolution in hyperbolic space.** The desingularisation construction that will be carried out in this paper will be controlled by two positive parameters. The first, which we will denote by $\epsilon$, will represent the factor by which the Costa-Hoffman-Meeks surface will be rescaled, and the second, which we will denote by $R$, will be such that the gluing operation will be performed over an annulus of inner and outer radii $\epsilon R$ and $2\epsilon R$ respectively. These parameters will be made to tend to zero and infinity respectively as follows. Given $\eta \ll 1$, we will suppose that

$$\epsilon R^{5-2\eta} \leq \frac{1}{\Lambda} \text{ and } \epsilon R^{5-\eta} \geq \Lambda, \tag{2.1}$$

where $\Lambda \gg 1$ will be chosen arbitrarily large. Theorem A will follow upon showing that, for sufficiently large $\Lambda$, the desingularisation construction indeed yields a minimal surface. We will also make use of catenoidal ends which vary over compact families. For this reason, a constant $C \gg 1$ will be fixed, and the parameter $c \in \mathbb{R}$ will be chosen such that

$$|c| \leq C. \tag{2.2}$$

We will see in Lemma 2.1.2, below, how $c$ parametrizes catenoidal ends.

Throughout this paper, we will use the following explicit parametrisations of two- and three-dimensional hyperbolic spaces. First, using polar coordinates, the metric $g := g^1$ is defined over $\mathbb{R}^2$ by

$$g := dr^2 + \sinh^2(r)d\theta^2, \tag{2.3}$$

where $r$ here denotes euclidean distance from the origin. The riemannian manifold $(\mathbb{R}^2, g)$ is then identified with two-dimensional hyperbolic space $\mathbb{H}^2$. Indeed, it is simply the parametrisation of $\mathbb{H}^2$ by geodesic coordinates about some point. In particular, all radial lines in this parametrisation are unit speed geodesics and the hyperbolic and euclidian distances to the origin coincide. The metric $\bar{g} := \bar{g}^1$ is defined over $\mathbb{R}^3$ by

$$\bar{g} := dr^2 + \sinh^2(r)d\theta^2 + \cosh^2(r)dt^2, \tag{2.4}$$

where $r$ now denotes euclidean distance from the $z$-axis. The riemannian manifold $(\mathbb{R}^3, \bar{g})$ is identified with three-dimensional hyperbolic space $\mathbb{H}^3$. Indeed, it is the parametrisation of $\mathbb{H}^3$ by Fermi coordinates about the geodesic $\Gamma := \{(0,0)\} \times \mathbb{R}$. In particular, the following properties stand out: every horizontal plane in $\mathbb{R}^3$ is totally geodesic; for every horizontal geodesic $\bar{\Gamma}$ passing through the origin, the cartesian product $\bar{\Gamma} \times \Gamma$ is totally geodesic; and, for all $R$, if $C(R)$ denotes the geodesic circle of radius $R$ about the origin in $\mathbb{R}^2$, then the surface $C(R) \times \mathbb{R}$ is a hyperbolic cylinder of radius $R$ about $\Gamma$.

Consider now a rotationally symmetric minimal surface $\Sigma$ in $\mathbb{H}^3$ which is a graph of some function $u$ over some annulus $A(a, \infty)$. We will henceforth refer to such surfaces as *hyperbolic minimal ends* in order to distinguish them from other types of minimal ends that will be introduced in Section 3.1, below.

**Lemma 2.1.1**

*The function $u$ satsfies*

$$u_r = \frac{(F/\pi)}{\cosh(r)\sqrt{\sinh^2(2r) - (F/\pi)^2}}, \tag{2.5}$$

where $F$ is a positive constant, which we henceforth refer to as the *flux* of $\Sigma$. Furthermore, $F$ satisfies

$$F = \pi \sinh(2r_0),$$

where $r_0$ is the neck radius of $\Sigma$.

**Remark:** The neck radius of $\Sigma$ can be defined to be the infimal value of $r$ such that $\Sigma$ extends to a smooth, rotationally symmetric graph over the annulus $A(r, \infty)$.

**Proof:** Indeed, for all $R$, let $\hat{C}(R)$ denote the intersection of $\Sigma$ with the cylinder $C(R) \times \mathbb{R}$. Given a Killing field $X$ of $\mathbb{H}^3$, the flux of $\Sigma$ around this circle is given by

$$F(R) = \int_{\hat{C}(R)} \langle X, \nu \rangle dl,$$

where $dl$ here denotes the length element of $\hat{C}(R)$, and $\nu$ denotes the unit conormal vector field along this curve. Letting $X := (0, 0, 1)$ be the Killing field of translations in the vertical direction, we obtain

$$F(R) = \frac{\pi \sinh(2R) \cosh(R) u_r(R)}{\sqrt{1 + \big(\cosh(R) u_r(R)\big)^2}}.$$

Since this quantity is constant whenever $\Sigma$ is minimal (c.f. [19]), the first result follows. Finally, the neck radius is defined to be the infimal distance to $\Gamma$ over $\Sigma$. Since this is realized when $u_r$ is infinite, the second result follows, and this completes the proof. $\square$

**Lemma 2.1.2**

Let $u : [\epsilon R, \epsilon R^4] \to \mathbb{R}$ satisfy (2.5) with flux $F = 2\pi \epsilon c$ for some $(\epsilon, R)$ and $c$ satisfying (2.1) and (2.2) respectively. Then, for sufficiently large $\Lambda$,

$$u_r = \frac{\epsilon c}{r} + O\left(\left(r + \frac{\epsilon}{r}\right)^3 \frac{1}{r^k}\right). \tag{2.6}$$

**Proof:** Indeed, by Taylor's Theorem,

$$\frac{1}{\cosh(r)} = 1 + O\big(r^{2-k}\big) = 1 + O\left(\left(r + \frac{\epsilon}{r}\right)^2 \frac{1}{r^k}\right).$$

Likewise,

$$\sinh^2(2r) = 4r^2 + O\big(r^{4-k}\big),$$

so that

$$\sinh^2(2r) - \left(\frac{F}{\pi}\right)^2 = 4r^2 \left[1 + O\left(\left(r + \frac{\epsilon}{r}\right)^2 \frac{1}{r^k}\right)\right].$$

Taking the reciprocal of the square root of this function yields

$$\frac{1}{\sqrt{\sinh^2(r) - \left(\frac{F}{\pi}\right)^2}} = \frac{1}{2r} \left[1 + O\left(\left(r + \frac{\epsilon}{r}\right)^2 \frac{1}{r^k}\right)\right],$$

and the result follows. $\square$

Differentiating (2.5) twice with respect to $c$ and repeating the argument of Lemma 2.1.2 likewise yields

**Lemma 2.1.3**

*For $(\epsilon, R)$ and $c$ satisfing (2.1) and (2.2) respectively, let $u_{c,r} : [\epsilon R, \epsilon R^4] \to \mathbb{R}$ be the unique solution of (2.5) with flux $F = 2\pi\epsilon c$. Then, for sufficiently large $\Lambda$,*

$$\frac{\partial u_{c,r}}{\partial c} = \frac{\epsilon}{r} + O\left(\left(r + \frac{\epsilon}{r}\right)^3 \frac{1}{r^k}\right), \text{ and}$$
$$\frac{\partial^2 u_{c,r}}{\partial c^2} = O\left(\left(r + \frac{\epsilon}{r}\right)^3 \frac{1}{r^k}\right). \tag{2.7}$$

**2.2 - The modified Jacobi operator.** Since the computations of this section will also be of use at later stages in the sequel, we will consider here the more general case of a graph $\Sigma$ of some function $u : A(\epsilon R, \epsilon R^4) \to \mathbb{R}$ which is not necessarily minimal nor rotationally symmetric, but which nonetheless satisfies

$$u = a + \epsilon c \text{Log}(r) + O\left(r\left(r + \frac{\epsilon}{r}\right)^3 \frac{1}{r^k}\right). \tag{2.8}$$

The Jacobi operator of $\Sigma$, which yields the infinitesimal variations of mean curvature that result from infinitesimal *normal* perturbations of the surface, is given, up to a sign, by (c.f. Appendix A.3)

$$Jf := \Delta^u f + \left[\text{Tr}\left((A^u)^2\right) + 2\text{Ric}(N^u, N^u)\right]f, \tag{2.9}$$

where Ric here denotes the Ricci curvature tensor* of $\mathbb{H}^3$ and $N^u$, $A^u$ and $\Delta^u$ denote respectively the upward-pointing unit normal vector field of the graph of $u$, its shape operator and its Laplace-Beltrami operator.

It turns out that the zero'th order coefficients of $J$ diverge rapidly as $\Lambda$ tends to infinity (as will become clearer in Section 2.5, below). For this reason, we choose to work with what we call the modified Jacobi operator. As we will see presently, the definition of this operator will vary slightly depending on the context. In the case at hand, it is defined as follows. First, let $\chi_1$ be the cut-off function of the transition region $A(1, 2)$. Define the vector field $X^u$ over $\Sigma$ by

$$X^u := \chi_1 e_z + (1 - \chi_1)N^u, \tag{2.10}$$

define the function $\psi^u : \Sigma \to \mathbb{R}$ by

$$\psi^u := \overline{g}\left(X^u, N^u\right) = \chi_1 \overline{g}(e_z, N^u) + (1 - \chi_1), \tag{2.11}$$

and, for $f \in C_0^\infty(\Sigma)$, define $\mathcal{E}_f : \Sigma \to \mathbb{H}^3$ by

$$\mathcal{E}_f(x) := x + tf(x)X^u(x).$$

---

* The Ricci curvature tensor is here normalised so that the unit sphere in Euclidean space has curvature equal to $\delta_{ij}$.

For all sufficiently small $f$, $\mathcal{E}_f$ is an embedding, and we define $\mathcal{H}_f : \Sigma \to \mathbb{R}$ to be its mean curvature function with respect to the metric $\overline{g}$ at the point $x$. The *modified Jacobi operator* of $\Sigma$ is then given by

$$(\hat{J}f)(x) := \frac{1}{\psi^u} \left. \frac{\partial}{\partial t} \mathcal{H}_{tf}(x) \right|_{t=0} . \tag{2.12}$$

Observe that this operator differs from the usual Jacobi operator merely by the inclusion of a factor of $1/\psi^u$ and the use of the vector field $X^u$ instead of the unit normal vector field in the definition of $\mathcal{E}$ (c.f. Appendix A.3). In particular, since $\Sigma$ is a graph, it is everywhere transverse to $X^u$, so that the function $\psi^u$ is everywhere strictly positive and the operator $\hat{J}$ is indeed well-defined.

The modified Jacobi operator is related to the usual Jacobi operator by

$$\hat{J}f = \left( \frac{1}{\psi^u} J \psi^u + \frac{1}{\psi^u} \pi^u(X^u) \mathcal{H}_0 \right) f, \tag{2.13}$$

where $\pi^u(X^u)$ here denotes the tangential component of $X^u$ and $\mathcal{H}_0$ denotes the mean curvature function of $\Sigma$. In particular, when $\Sigma$ is minimal, $\mathcal{H}_0$ vanishes, so that

$$\hat{J}f = \frac{1}{\psi^u} J(\psi^u f). \tag{2.14}$$

On the other hand, since $e_z$ is a Killing field of $\mathbb{H}^3$, the quantity

$$\hat{J} \cdot 1 = \left( \frac{1}{\psi^u} J \psi^u + \frac{1}{\psi^u} \pi^u(X^u) \mathcal{H}_0 \right)$$

vanishes over $\Sigma \cap (B(1) \times \mathbb{R})$, so that, by (A.7), over this region

$$\hat{J}f = \Delta^u f + \frac{2}{\mu^u} g^{u,ij} \mu_i^u f_j, \tag{2.15}$$

where the function $\mu^u$ is defined as in Appendix A.4. The relations (2.14) and (2.15) will be used at various points in what follows.

Since $\Sigma$ is a graph over $A(\epsilon R, 2\epsilon R^4)$, $\hat{J}$ may also be thought of as an operator acting on functions defined over this domain. Using this identification, we obtain

**Lemma 2.2.1**

*Over $A(\epsilon R, \epsilon R^4)$,*

$$\hat{J}f = \Delta^g f - \frac{\epsilon^2 c^2}{r^2} f_{rr} + 2\left( r + \frac{\epsilon^2 c^2}{r^3} \right) f_r + \mathcal{R}f, \tag{2.16}$$

*where $\Delta^g$ here denotes the Laplace-Beltrami operator of $\mathbb{H}^2$ and the remainder $\mathcal{R}$ is a second order linear partial differential operator given by*

$$\mathcal{R}f = a^{ij} f_{ij} + b^i f_i,$$

where

$$a = O\bigg(\bigg(r + \frac{\epsilon}{r}\bigg)^4 \frac{1}{r^k}\bigg), \text{ and}$$

$$b = O\bigg(\bigg(r + \frac{\epsilon}{r}\bigg)^4 \frac{1}{r^{k+1}}\bigg). \tag{2.17}$$

**Proof:** First, combining (2.8) and (A.8), yields

$$\mu^u = 1 + \frac{r^2}{2} - \frac{\epsilon^2 c^2}{2r^2} + O\bigg(\bigg(r + \frac{\epsilon}{r}\bigg)^4 \frac{1}{r^k}\bigg).$$

Next, given a point $x$ in $\mathbb{H}^2$, consider the orthonormal basis $(e_r, e_\theta)$, where $e_r$ is the unit vector in the radial direction pointing away from the origin and $e_\theta$ is the hyperbolic unit vector in the counterclockwise angular direction around the origin. By (2.8) again and (A.11), with respect to this basis,

$$g^u(e_i, e_j) = \begin{pmatrix} 1 & \\ & 1 \end{pmatrix} + \begin{pmatrix} \epsilon^2 c^2/r^2 & \\ & 0 \end{pmatrix} + O\bigg(\bigg(r + \frac{\epsilon}{r}\bigg)^4 \frac{1}{r^k}\bigg).$$

Likewise, the Hessian of $u$ with respect to $g$ is given in this basis by,

$$\text{Hess}^g(u)(e_i, e_j) = \begin{pmatrix} -\epsilon c/r^2 & \\ & \epsilon c/r^2 \end{pmatrix} + O\bigg(\bigg(r + \frac{\epsilon}{r}\bigg)^3 \frac{1}{r^{k+1}}\bigg).$$

In particular, since the leading order term in the Hessian has zero trace, combining the previous two formulae yields

$$\cosh^2(r) g^{u,ij} g^{u,kp} \text{Hess}^g(u)_{ij} u_p = O\bigg(\bigg(r + \frac{\epsilon}{r}\bigg)^4 \frac{1}{r^{k+1}}\bigg).$$

Likewise,

$$\cosh(r)\sinh(r) g^{u,ij} g^{u,kp}\big(u_i u_p r_j + u_j u_p r_i - u_i u_j r_p\big) = O\bigg(\bigg(r + \frac{\epsilon}{r}\bigg)^4 \frac{1}{r^{k+1}}\bigg),$$

so that, by (A.15), if $f$ is a smooth function defined over the domain of $u$, then its Laplacian with respect to $g^u$ satisfies

$$\Delta^u f = g^{u,ij} \text{Hess}^g(f)_{ij} + \mathcal{R}_1 f,$$

$$= \Delta^g f - \frac{\epsilon^2 c^2}{r^2} f_{rr} + \mathcal{R}_2 f,$$

where the remainders $\mathcal{R}_1$ and $\mathcal{R}_2$ are second order linear partial differential operators of the desired form. Finally, we have

$$\frac{2}{\mu^u} g^{u,ij} \mu_i = 2\bigg(r + \frac{\epsilon^2 c^2}{r^3}\bigg)\delta_r^j + O\bigg(\bigg(r + \frac{\epsilon}{r}\bigg)^4 \frac{1}{r^{k+1}}\bigg),$$

and the result now follows by (2.15). $\square$

**2.3 - The Sobolev, Hölder and hybrid norms.** Let $m$ be a non-negative integer and $\alpha \in [0, 1]$ a real number. Let $\| \cdot \|_{H^m(\mathbb{H}^2)}$ and $\| \cdot \|_{C^{m,\alpha}(\mathbb{H}^2)}$ denote respectively the *Sobolev norm* of order $m$ and the *Hölder norm* of order $(m, \alpha)$ of functions over $\mathbb{H}^2$, as defined in Appendix A.5. Let $H^m(\mathbb{H}^2)$ denote the *Sobolev space* of measurable functions $f$ whose distributional derivatives up to and including order $m$ are square integrable over $\mathbb{H}^2$ and let $C^{m,\alpha}(\mathbb{H}^2)$ denote the *Hölder space* of $m$-times differentiable functions $f$ over $\mathbb{H}^2$ which satisfy $\|f\|_{C^{m,\alpha}(\mathbb{H}^2)} < \infty$.

We modify these classical function spaces by a continuous family of weights as follows. Let $\chi_2$ denote the cut-off function of the transition region $A(2, 4)$ and, for $\gamma \in \mathbb{R}$, define the weight $w_\gamma : \mathbb{R}^2 \to \mathbb{R}$ by

$$w_\gamma := \chi_2 + (1 - \chi_2) e^{\gamma r}, \tag{2.18}$$

where $r$ here denotes the distance from the origin. With $m$ and $\alpha$ as before, the *weighted Sobolev and Hölder norms* with weight $\gamma$ are defined by

$$\begin{aligned} \|\phi\|_{H^m_\gamma(\mathbb{H}^2)} &:= \|w_\gamma \phi\|_{H^m(\mathbb{H}^2)}, \text{ and} \\ \|\phi\|_{C^{m,\alpha}_\gamma(\mathbb{H}^2)} &:= \|w_\gamma \phi\|_{C^{m,\alpha}(\mathbb{H}^2)}, \end{aligned} \tag{2.19}$$

and the *weighted Sobolev and Hölder spaces* with weight $\gamma$ are defined by

$$\begin{aligned} H^m_\gamma(\mathbb{H}^2) &:= \left\{ w^{-1}_\gamma f \mid f \in H^m(\mathbb{H}^2) \right\}, \text{ and} \\ C^{m,\alpha}_\gamma(\mathbb{H}^2) &:= \left\{ w^{-1}_\gamma f \mid f \in C^{m,\alpha}(\mathbb{H}^2) \right\}. \end{aligned} \tag{2.20}$$

These spaces are trivially Banach spaces and multiplication by $w_\gamma$ defines linear isomorphisms from $H^m_\gamma(\mathbb{H}^2)$ into $H^m(\mathbb{H}^2)$ and from $C^{m,\alpha}_\gamma(\mathbb{H}^2)$ into $C^{m,\alpha}(\mathbb{H}^2)$.

Consider now a complete, horizontal, totally geodesic plane $\Sigma$ in $\mathbb{H}^3$. Its Jacobi operator is given by

$$Jf := J_0 f := \Delta^g f - 2f, \tag{2.21}$$

where $\Delta^g$ here denotes the Laplace-Beltrami operator of $\mathbb{H}^2$. By minimality and (2.14), its modified Jacobi operator is given by

$$\hat{J}_0 f := \frac{1}{\psi^u} J_0(\psi^u f),$$

where $\psi^u$ is the function given by (2.11). Likewise, by (2.15), over $B(1)$, its modified Jacobi operator is also given by

$$\begin{aligned} \hat{J}_0 f &= \Delta^g f + 2\tanh(r) f_r \\ &= \Delta^g f + 2r f_r + \mathcal{R} f \end{aligned}, \tag{2.22}$$

where the remainder $\mathcal{R}$ is a second-order linear partial differential operator of the form

$$\mathcal{R}^{ij} f = a^{ij} f_{ij} + b^i f_i,$$

and whose coefficients satisfy (2.17).

**Theorem 2.3.1**

*(1) For all sufficiently small $\gamma$, $\hat{J}_0$ defines a linear isomorphism from $H^2_\gamma(\mathbb{H}^2)$ into $H^0_\gamma(\mathbb{H}^2)$; and*

*(2) for all $\alpha \in ]0, 1[$, and for all sufficiently small $\gamma$, $\hat{J}_0$ defines a linear isomorphism from $C^{2,\alpha}_\gamma(\mathbb{H}^2)$ into $C^{0,\alpha}_\gamma(\mathbb{H}^2)$.*

**Proof:** We only prove (1), as the proof of (2) is identical. Since $\psi^u$ is strictly positive, smooth and equal to 1 outside of some compact set, multiplication by $\psi^u$ defines linear isomorphisms from $H^2_\gamma(\mathbb{H}^2)$ to itself and from $H^0_\gamma(\mathbb{H}^2)$ to itself. It thus suffices to prove the result for $J_0$. Next, observe that $J_0$ defines a linear isomorphism from $H^2_\gamma(\mathbb{H}^2)$ into $H^0_\gamma(\mathbb{H}^2)$ if and only if $J_{0,\gamma} := M_\gamma J_0 M_\gamma^{-1}$ defines a linear isomorphism from $H^2(\mathbb{H}^2)$ into $H^0(\mathbb{H}^2)$, where $M_\gamma$ here denotes the operator of multiplication by $w_\gamma$. Since this family varies continuously in the operator norm, it suffices to prove the result for $\gamma = 0$. However, this final case follows by the classical theory of elliptic operators (c.f. [5]), and this completes the proof. $\square$

Many of the estimates derived in the sequel are best expressed in terms of a hybrid norm which we now introduce. For all non-negative integer $m$, for all $\alpha \in [0, 1]$ and for all real $\gamma$, we define $\| \cdot \|_{m,\alpha,\gamma}$, the *hybrid norm* with weight $\gamma$ of functions over $\mathbb{H}^2$ by

$$\|f\|_{m,\alpha,\gamma} := \|f\|_{C^{m,\alpha}_\gamma(\mathbb{H}^2)} + \frac{1}{(\epsilon R)} \|f\|_{H^m_\gamma(\mathbb{H}^2)}. \tag{2.23}$$

We denote by $\mathcal{L}^{m,\alpha}_\gamma(\mathbb{H}^2)$ the Banach space of $m$-times differentiable functions $f$ over $\mathbb{H}^2$ which satisfy $\|f\|_{m,\alpha,\gamma} < \infty$. In terms of the hybrid norm, Theorem 2.3.1 immediately yields

**Theorem 2.3.2**

*For all $\alpha \in ]0, 1[$ and for all sufficiently small $\gamma$, $\hat{J}_0$ defines a linear isomorphism from $\mathcal{L}^{2,\alpha}_\gamma(\mathbb{H}^2)$ into $\mathcal{L}^{0,\alpha}_\gamma(\mathbb{H}^2)$. Furthermore, the operator norms of $\hat{J}_0$ and its inverse are uniformly bounded independent of $\Lambda$.*

We now extend Theorem 2.3.2 to the case where $\Sigma$ is a hyperbolic minimal end with flux $F = 2\pi\epsilon c$ which is a graph of some function $u$ over the annulus $A(\epsilon R, \infty)$. Recall first that the modified Jacobi operator $\hat{J}$ of $\Sigma$ may be considered as an operator acting on functions over $A(\epsilon R, \infty)$. However, it will now be useful to extend it to one defined over the whole of $\mathbb{R}^2$ as follows. First, given a function $f : A(\epsilon R, \infty) \to \mathbb{R}$, its *canonical extension* $\tilde{f} : \mathbb{R}^2 \to \mathbb{R}$ is defined such that

(1) $\tilde{f}(x) := f(x)$ over $A(\epsilon R, \infty)$;

(2) $\tilde{f}(0)$ is equal to the mean value of $f$ over the circle $C(\epsilon R)$; and

(3) $\tilde{f}$ restricts to a linear function over every radial line in $B(\epsilon R)$.

Given a linear operator $L$ over $A(\epsilon R, \infty)$, its *canonical extension* $\tilde{L}$ is defined to be the operator over $\mathbb{R}^2$ whose coefficients are the canonical extensions of each of the coefficients of $L$. Observe that if a function $f$ is Lipschitz, then so too is its canonical extension and

$$\|\tilde{f}\|_{C^{0,1}} \lesssim \|f\|_{C^{0,1}}.$$

Likewise, if an operator $L$ has any rotational symmetries, then so too does its canonical extension. We will henceforth identify all operators with their canonical extensions, and a perturbation of Theorem 2.3.2 will then yield

**Theorem 2.3.3**

*For all sufficiently small $\alpha \in ]0, 1[$, for all sufficiently small $\gamma$ and for all sufficiently large $\Lambda$, the modified Jacobi operator $\hat{J}$ defines a linear isomorphism from $\mathcal{L}_\gamma^{2,\alpha}(\mathbb{H}^2)$ into $\mathcal{L}_\gamma^{0,\alpha}(\mathbb{H}^2)$. Furthermore, the operator norms of $\hat{J}$ and its inverse are uniformly bounded as $\Lambda$ tends to infinity.*

This result follows immediately from Theorem 2.3.2 and Lemmas 2.4.1 and 2.5.1 below.

**2.4 - The regular component.** Theorem 2.3.3 will be proven by showing that $\hat{J}$ converges towards $\hat{J}_0$ in every operator norm of interest to us as $\Lambda$ tends to infinity. To this end, we decompose the difference $(\hat{J} - \hat{J}_0)$ into two components which we analyse separately. Thus, let $\chi_{\epsilon R^4}$ be the cut-off function of the transition region $A(\epsilon R^4, 2\epsilon R^4)$, and define

$$\mathcal{M}f := 2\chi_{\epsilon R^4} \frac{\epsilon^2 c^2}{r^3} f_r, \text{ and}$$
$$\mathcal{N}f := \left( (\hat{J} - \hat{J}_0) - \mathcal{M} \right) f. \tag{2.24}$$

We refer to $\mathcal{M}$ and $\mathcal{N}$ respectively as the *singular component* and the *regular component* of the difference $(\hat{J} - \hat{J}_0)$.

**Lemma 2.4.1**

*(1) For sufficiently small $\gamma$, the operator norm of $\mathcal{N}$ considered as a linear map from $H_\gamma^2(\mathbb{H}^2)$ into $H_\gamma^0(\mathbb{H}^2)$ tends to zero as $\Lambda$ tends to infinity; and*

*(2) for sufficiently small $\alpha \in [0, 1]$, and for sufficiently small $\gamma$, the operator norm of $\mathcal{N}$, considered as a linear map from $C_\gamma^{2,\alpha}(\mathbb{H}^2)$ into $C_\gamma^{0,\alpha}(\mathbb{H}^2)$ tends to zero as $\Lambda$ tends to infinity.*

Denote

$$\mathcal{N}f =: a^{ij} f_{ij} + b^i f_i. \tag{2.25}$$

Lemma 2.4.1 will follow from Lemmas 2.4.2 and 2.4.3, below.

**Lemma 2.4.2**

*The coefficients $a$ and $b$ satisfy*

$$\left\| a|_{A(2\epsilon R^4, \infty)} \right\|_{C^{0,1}(\mathbb{H}^2)}, \left\| b|_{A(2\epsilon R^4, \infty)} \right\|_{C^{0,1}(\mathbb{H}^2)} \to 0$$

*as $\Lambda$ tends to infinity.*

**Proof:** Denote $f = \cosh(r)u_r$. By (2.5),

$$f = \frac{(F/\pi)}{\sqrt{\sinh^2(2r) - (F/\pi)^2}}.$$

Differentiating twice yields

$$f_r = -(\pi/F)^2 \sinh(4r)f^3, \text{ and}$$
$$f_{rr} = -4(\pi/F)^2 f^3 - 8(\pi/F)^2 \sinh^2(2r)f^3 + 3(\pi/F)^4 \sinh^2(4r)f^5.$$

We now study these functions over the annulus $A(2\epsilon R^4, \infty)$. Since $F = 2\pi\epsilon c$, for $\Lambda$ sufficiently large,

$$\sqrt{\sinh^2(2r) - (F/\pi)^2} \gtrsim \sinh(2r),$$

so that, over this annulus, for $k \in \{0, 1, 2\}$,

$$\left| \frac{\partial^k f}{\partial r^k} \right| \lesssim \frac{\epsilon}{\sinh^{k+1}(2r)} + \frac{\epsilon}{\sinh(2r)}.$$

Since

$$u_r = \frac{1}{\cosh(r)} f,$$

it then follows by rotational symmetry that, over $A(2\epsilon R^4, \infty)$, for $k \in \{1, 2, 3\}$,

$$\left\| (\nabla^g)^k u \right\|_g \lesssim \frac{\epsilon}{\sinh^k(2r)} + \frac{\epsilon}{\sinh(2r)} \lesssim \frac{1}{\epsilon^{k-1}R^{4k}}. \tag{2.26}$$

The result now follows since, by (2.9), (2.10), (2.11) and (2.12), the coefficients of $\hat{J}$ only depend on the first and second derivatives of $u$. $\square$

**Lemma 2.4.3**

*For sufficiently small $\alpha \in [0, 1]$, the coefficients $a$ and $b$ satisfy*

$$\left\| a|_{B(2\epsilon R^4)} \right\|_{C^{0,\alpha}}, \left\| b|_{B(2\epsilon R^4)} \right\|_{C^{0,\alpha}} \to 0$$

*as $\Lambda$ tends to infinity.*

**Remark:** Observe that this lemma fails when the exponent of $R$ in the second inequality in (2.1) is greater than or equal to 5.

**Proof:** Indeed, by (2.16), (2.22) and (2.24), over the annulus $A(\epsilon R, 2\epsilon R^4)$,

$$a^{ij} = -\frac{\epsilon^2 c^2 x^i x^j}{r^4} + \mathrm{O}\left( \left( r + \frac{\epsilon}{r} \right)^4 \frac{1}{r^k} \right), \text{ and}$$
$$b^i = 2\left(1 - \chi_{\epsilon R^4}\right)\frac{\epsilon^2 c^2 x^i}{r^4} + \mathrm{O}\left( \left( r + \frac{\epsilon}{r} \right)^4 \frac{1}{r^{k+1}} \right).$$

It follows that

$$|a| \lesssim \frac{\epsilon^2}{r^2} + r^4, \text{ and}$$
$$|Da| \lesssim \frac{\epsilon^2}{r^3} + r^3,$$

so that, by (2.1),

$$\left\| a|_{A(\epsilon R, 2\epsilon R^4)} \right\|_{C^0} \lesssim \frac{1}{R^2} + \epsilon^4 R^{16} \lesssim \frac{1}{R^2}, \text{ and}$$

$$\left[ a|_{A(\epsilon R, 2\epsilon R^4)} \right]_1 \lesssim \frac{1}{\epsilon R^3} + \epsilon^3 R^{12} \lesssim \frac{1}{\epsilon R^3}.$$

Since $a$ is extended canonically over $B(\epsilon R)$, these inequalities in fact hold over the whole of $B(2\epsilon R^4)$ and so, by (A.23),

$$\left[ a|_{B(2\epsilon R^4)} \right]_{\frac{1}{2}} \leq \sqrt{2} \left\| a|_{B(2\epsilon R^4)} \right\|_{C^0}^{\frac{1}{2}} \left[ a|_{B(2\epsilon R^4)} \right]_1^{\frac{1}{2}} \lesssim \frac{1}{\sqrt{\epsilon R^5}}.$$

Bearing in mind that $(1 - \chi_{\epsilon R^4}) = \mathrm{O}(r^{-k})$, over the annulus $A(\epsilon R^4, 2\epsilon R^4)$, we obtain

$$|b| \lesssim \frac{\epsilon^2}{r^3} + r^3, \text{ and}$$

$$|Db| \lesssim \frac{\epsilon^2}{r^4} + r^2,$$

so that, by (2.1),

$$\left\| b|_{A(\epsilon R^4, 2\epsilon R^4)} \right\|_{C^0} \lesssim \frac{1}{\epsilon R^{12}} + \epsilon^3 R^{12} \lesssim \frac{1}{\epsilon R^5}, \text{ and}$$

$$\left[ b|_{A(\epsilon R^4, 2\epsilon R^4)} \right]_1 \lesssim \frac{1}{\epsilon^2 R^{16}} + \epsilon^2 R^8 \lesssim \frac{1}{(\epsilon R^3)^2}.$$

On the other hand, over the annulus $A(\epsilon R, \epsilon R^4)$,

$$|b| \lesssim r^3 + \frac{\epsilon^4}{r^5}, \text{ and}$$

$$|Db| \lesssim r^2 + \frac{\epsilon^4}{r^6},$$

so that, by (2.1),

$$\left\| b|_{A(\epsilon R, \epsilon R^4)} \right\|_{C^0} \lesssim \epsilon^3 R^{12} + \frac{1}{\epsilon R^5} \lesssim \frac{1}{\epsilon R^5}, \text{ and}$$

$$\left[ b|_{A(\epsilon R, \epsilon R^4)} \right]_1 \lesssim \epsilon^2 R^8 + \frac{1}{\epsilon^2 R^6} \lesssim \frac{1}{(\epsilon R^3)^2}.$$

As before, since $b$ is extended canonically over $B(\epsilon R)$, these inequalities continue to hold over the whole of $B(2\epsilon R^4)$. Thus, by (A.23), for for $\alpha \in [0,1]$,

$$\left[ b|_{B(2\epsilon R^4)} \right]_\alpha \leq 2^{1-\alpha} \left\| b|_{B(2\epsilon R^4)} \right\|_{C^0}^{1-\alpha} \left[ b|_{B(2\epsilon R^4)} \right]_1^\alpha \lesssim \frac{1}{\epsilon^{1+\alpha} R^{5+\alpha}}.$$

Choosing $\alpha$ such that $(5 + \alpha)/(1 + \alpha) \geq 5 - \eta$, the result follows for $b$. $\square$

## 2.5 - The singular component.

We now turn our attention to the operator norm of $\mathcal{M}$.

**Lemma 2.5.1**

*For sufficiently small $\alpha \in [0,1]$ and for sufficiently small $\gamma$, the operator norm of $\mathcal{M}$ considered as a linear map from $\mathcal{L}^{2,\alpha}_\gamma$ into $\mathcal{L}^{0,\alpha}_\gamma$ tends to 0 as $\Lambda$ tends to infinity.*

In contrast to the regular component, the coefficients of $\mathcal{M}$ do not tend pointwise to zero as $\Lambda$ tends to infinity. However, we will show that this is compensated for by the fact that $\mathcal{M}$ is supported over a ball of small diameter. Indeed, by definition of the canonical extension,

$$\mathcal{M}\phi =: a^i \phi_i, \tag{2.27}$$

where

$$a^i = \begin{cases} 2\chi_{\epsilon R^4} \dfrac{\epsilon^2 c^2 x_i}{r^4} & \text{if } \|x\| > \epsilon R, \text{ and} \\[2mm] \dfrac{2c^2 x_i}{\epsilon^2 R^4} & \text{if } \|x\| \leq \epsilon R. \end{cases} \tag{2.28}$$

The required estimates will follow from Lemmas 2.5.3 and 2.5.4, below. These lemmas both rely on the following key estimate, which reveals the significance of the hybrid norm.

**Lemma 2.5.2**

*For sufficiently small $\alpha$ and for sufficiently small $\gamma$,*

$$\|f\|_{C^{1,\alpha}_\gamma(\mathbb{H}^2)} \lesssim (\epsilon R)^{1-2\alpha} \|f\|_{2,\alpha,\gamma}. \tag{2.29}$$

**Remark:** It will be useful in the sequel to observe that this relation also holds for spaces of functions defined over an unbounded annulus.

**Proof:** Indeed, by the Sobolev embedding theorem, for all $\beta < 1$,

$$\|f\|_{C^{0,\beta}_\gamma(\mathbb{H}^2)} \lesssim \|f\|_{H^2_\gamma(\mathbb{H}^2)} \lesssim (\epsilon R)\|f\|_{2,\alpha,\gamma}.$$

Setting $\beta = (1-\alpha)$ and using (A.23) and (A.24), we obtain

$$\|f\|_{C^{1,\alpha}_\gamma(\mathbb{H}^2)} \lesssim (\epsilon R)^{\frac{1}{1+2\alpha}} \|f\|_{2,\alpha,\gamma} \lesssim (\epsilon R)^{1-2\alpha} \|f\|_{2,\alpha,\gamma},$$

as desired. $\square$

**Lemma 2.5.3**

*For sufficiently small $\alpha \in [0,1]$ and for sufficiently small $\gamma$, the operator norm of $\mathcal{M}$, considered as a map from $\mathcal{L}^{2,\alpha}_\gamma(\mathbb{H}^2)$ into $C^{0,\alpha}_\gamma(\mathbb{H}^2)$ tends to 0 as $\Lambda$ tends to infinity.*

**Proof:** Indeed, over $A(\epsilon R, 2\epsilon R^4)$,

$$a^i = \mathrm{O}\left(\frac{\epsilon^2}{r^{3+k}}\right),$$

so that

$$\big\| a^i|_{A(\epsilon R, 2\epsilon R^4)} \big\|_{C^0} \lesssim \frac{1}{\epsilon R^3}, \text{ and}$$

$$\big[ a^i|_{A(\epsilon R, 2\epsilon R^4)} \big]_1 \lesssim \frac{1}{\epsilon^2 R^4}.$$

Since $a^i$ is extended canonically over $B(\epsilon R)$, these inequalities in fact hold over the whole of $B(2\epsilon R^4)$ so that, by (A.23), for all $\alpha \in [0, 1]$,

$$[a^i]_\alpha \lesssim \frac{1}{(\epsilon R)^\alpha \epsilon R^3}.$$

It follows by (2.29) and (A.25) that

$$\|\mathcal{M}f\|_{C^{0,\alpha}_\gamma(\mathbb{H}^2)} \lesssim \frac{1}{(\epsilon R)^\alpha \epsilon R^3} \|f\|_{C^{1,\alpha}_\gamma(\mathbb{H}^2)} \lesssim \frac{1}{(\epsilon R)^{3\alpha} R^2} \|f\|_{2,\alpha,\gamma},$$

and the result follows by (2.1). $\square$

**Lemma 2.5.4**

*For sufficiently small $\alpha \in [0, 1]$ and for sufficiently small $\gamma$, the operator norm of $(\epsilon R)^{-1}\mathcal{M}$ considered as a map from $\mathcal{L}^{2\alpha}_\gamma(\mathbb{H}^2)$ into $H^0_\gamma(\mathbb{H}^2)$ tends to 0 as $\Lambda$ tends to infinity.*

**Proof:** Indeed, a direct calculation yields

$$\|a^i\|_{L^2_\gamma(\mathbb{H}^2)} \lesssim \frac{1}{R^2}.$$

Thus, bearing in mind (2.29),

$$\|(\epsilon R)^{-1}\mathcal{M}f\|_{L^2_\gamma(\mathbb{H}^2)} \lesssim (\epsilon R)^{-1}\|a^i\|_{L^2_\gamma(\mathbb{H}^2)} \|Df\|_{L^\infty}$$

$$\lesssim (\epsilon R)^{-1}\|a^i\|_{L^2_\gamma(\mathbb{H}^2)} \|f\|_{C^{1,\alpha}_\gamma(\mathbb{H}^2)}$$

$$\lesssim \frac{1}{(\epsilon R)^{2\alpha} R^2} \|f\|_{2,\alpha,\gamma},$$

and the result follows by (2.1). $\square$

**2.6 - Second order variations.** Consider again the functional $\mathcal{H}_f$ introduced at the start of Section 2.2. In addition to the above estimates, the proof of Theorem A will also require estimates for the second-order variation of this functional. We show

**Theorem 2.6.1**

*For sufficiently small $\alpha \in [0,1]$, for sufficiently small $\gamma$ and for sufficiently large $\Lambda$, if*

$$\big\| f|_{A(\epsilon R, \infty)} \big\|_{C^{1,\alpha}_\gamma(\mathbb{H}^2)} \leq \frac{1}{(\epsilon R)^\alpha R}, \tag{2.30}$$

*then*

$$\big\| \big(\mathcal{H}_f - \mathcal{H}_0 - \hat{J}f\big)|_{A(\epsilon R, \infty)} \big\|_{0,\alpha,\gamma} \lesssim (\epsilon R)^{1-2\alpha} \big\| f|_{A(\epsilon R, \infty)} \big\|^2_{2,\alpha,\gamma}. \tag{2.31}$$

**Remark:** Although we expect a finer analysis over the region $A(\epsilon R^4, \infty)$ to yield a stronger estimate, the above result is quite sufficient for our purposes.

**Proof:** Denote

$$\mathcal{R} := \mathcal{H}_f - \mathcal{H}_0 - \hat{J}f,$$

and consider first the restriction of this function to the annulus $A(\epsilon R^4, \infty)$. Since $\mathcal{H}$ is a second-order quasi-linear operator whose coefficients are uniformly bounded over this region independent of $\Lambda$, we have by Taylor's Theorem and (2.29) that

$$\big\| \mathcal{R}|_{A(\epsilon R^4, \infty)} \big\|_{0,\alpha,\gamma} \lesssim \big\| f|_{A(\epsilon R, \infty)} \big\|_{C^{1,\alpha}_\gamma} \big\| f|_{A(\epsilon R, \infty)} \big\|_{2,\alpha,\gamma}$$

$$\lesssim (\epsilon R)^{1-2\alpha} \big\| f|_{A(\epsilon R, \infty)} \big\|^2_{2,\alpha,\gamma}.$$

Consider now the restriction of $\mathcal{R}$ to the annulus $A(\epsilon R, \epsilon R^4)$. As in the previous sections, we now require a more careful analysis in order to account for singular behaviour as $\Lambda$ tends to infinity. To this end, consider two functions $\phi, \psi \in \mathcal{L}^{2,\alpha}_\gamma(\mathbb{H}^2)$ and suppose in addition that $\phi$ satisfies (2.30). Let $\hat{J}_\phi$ be the modified Jacobi operator of the graph of $(u + \phi)$. Observe, in particular, that $\hat{J}_0 = \hat{J}$. It suffices to show that, for all $\phi$ and for all $\psi$,

$$\big\| \big(\hat{J}_\phi \psi - \hat{J}\psi\big)|_{A(\epsilon R, \epsilon R^4)} \big\|_{0,\alpha,\gamma} \lesssim \frac{\epsilon}{(\epsilon R)^{5\alpha}} \|\phi\|_{2,\alpha,\gamma} \|\psi\|_{2,\alpha,\gamma}. \tag{2.32}$$

Indeed, substituting $\phi = tf$ and $\psi = f$ and integrating (2.32) with respect to $t$ over the interval $[0,1]$ yields

$$\big\| \mathcal{R}|_{A(\epsilon R, \epsilon R^4)} \big\|_{0,\alpha,\gamma} \lesssim \frac{\epsilon}{(\epsilon R)^{5\alpha}} \|f\|^2_{2,\alpha,\gamma},$$

and the result follows.

We now prove (2.32). By (2.15) and the formulae in Appendix A.4, over $B(1)$, the operator $\hat{J}_\phi$ is given by

$$\hat{J}_\phi \psi = F(x, D(u + \phi) \otimes D(u + \phi))^{ij} \text{Hess}^g(\psi)_{ij} + G(x, D(u + \phi) \otimes D(u + \phi))^i \psi_i$$
$$+ H(x, D(u + \phi))^{ijk} \text{Hess}^g(u + \phi)_{ij} \psi_k,$$

where $F$, $G$ and $H$ are smooth functions of their arguments and

$$H(x, 0)^{ijk} = 0.$$

It follows that

$$\hat{J}_\phi\psi - \hat{J}\psi = a^{ij}\text{Hess}^g(\psi)_{ij} + b^i\psi_i + c^{ijk}\text{Hess}^g(\phi)_{ij}\psi_k + d^{ijk}\text{Hess}^g(u)_{ij}\psi_k,$$

where

$$a^{ij} := F(x, D(u+\phi) \otimes D(u+\phi))^{ij} - F(x, Du \otimes Du)^{ij},$$
$$b^i := G(x, D(u+\phi) \otimes D(u+\phi))^i - G(x, Du \otimes Du)^i,$$
$$c^{ijk} := H(x, D(u+\phi))^{ijk}, \text{ and}$$
$$d^{ijk} := H(x, D(u+\phi))^{ijk} - H(x, Du)^{ijk}.$$

Bearing in mind (2.1), (2.8), (2.29), (2.30) and (A.25), we obtain

$$\left\|\left(D(u+\phi) \times D(u+\phi) - Du \otimes Du\right)|_{A(\epsilon R, \epsilon R^4)}\right\|_{C_\gamma^{0,\alpha}}$$
$$\lesssim \left(\left\|Du|_{A(\epsilon R, \epsilon R^4)}\right\|_{C_\gamma^{0,\alpha}} + \|D\phi\|_{C_\gamma^{0,\alpha}}\right)\|D\phi\|_{C_\gamma^{0,\alpha}}$$
$$\lesssim \frac{1}{(\epsilon R)^\alpha R}\|D\phi\|_{C_\gamma^{0,\alpha}}$$
$$\lesssim \frac{\epsilon}{(\epsilon R)^{3\alpha}}\|\phi\|_{2,\alpha,\gamma},$$

so that, by (A.27),

$$\left\|a^{ij}|_{A(\epsilon R, \epsilon R^4)}\right\|_{C_\gamma^{0,\alpha}}, \left\|b^i|_{A(\epsilon R, \epsilon R^4)}\right\|_{C_\gamma^{0,\alpha}} \lesssim \frac{\epsilon}{(\epsilon R)^{3\alpha}}\|\phi\|_{2,\alpha,\gamma}.$$

Likewise, using in addition the fact that $H(x,0)^{ijk} = 0$, we obtain

$$\left\|c^{ijk}|_{A(\epsilon R, \epsilon R^4)}\right\|_{C_\gamma^{0,\alpha}} \lesssim \frac{1}{(\epsilon R)^\alpha R} + \|D\phi\|_{C_\gamma^{0,\alpha}}$$
$$\lesssim \frac{1}{(\epsilon R)^\alpha R}, \text{ and}$$
$$\left\|d^{ijk}|_{A(\epsilon R, \epsilon R^4)}\right\|_{C_\gamma^{0,\alpha}} \lesssim (\epsilon R)^{1-2\alpha}\|\phi\|_{2,\alpha,\gamma}.$$

Finally, by (2.8) again

$$\left\|\text{Hess}^g(u)_{ij}|_{A(\epsilon R, \epsilon R^4)}\right\|_{C_\gamma^{0,\alpha}} \lesssim \frac{1}{(\epsilon R)^\alpha \epsilon R^2}, \text{ and}$$
$$\left\|\text{Hess}^g(u)_{ij}|_{A(\epsilon R, \epsilon R^4)}\right\|_{H_\gamma^0} \lesssim \frac{1}{R}.$$

Combining these relations yields

$$\left\|(\hat{J}_\phi\psi - \hat{J}\psi)|_{A(\epsilon R, \epsilon R^4)}\right\|_{C_\gamma^{0,\alpha}} \lesssim \frac{\epsilon}{(\epsilon R)^{5\alpha}}\|\phi\|_{2,\alpha,\gamma}\|\psi\|_{2,\alpha,\gamma}, \text{ and}$$
$$\frac{1}{(\epsilon R)}\left\|(\hat{J}_\phi\psi - \hat{J}\psi)|_{A(\epsilon R, \epsilon R^4)}\right\|_{H_\gamma^0} \lesssim \frac{\epsilon}{(\epsilon R)^{4\alpha}}\|\phi\|_{2,\alpha,\gamma}\|\psi\|_{2,\alpha,\gamma},$$

and the result follows. $\square$

## 3 - Surgery and the perturbation family.

**3.1 - The surgery operation.** We now rescale the parametrisations of two- and three-dimensional hyperbolic spaces introduced in Section 2.1. Thus, let $g^0$ and $\overline{g}^0$ be the euclidian metrics of $\mathbb{R}^2$ and $\mathbb{R}^3$ respectively and, for $\epsilon > 0$, let $g^\epsilon$ and $\overline{g}^\epsilon$ be the metrics defined respectively over these spaces by

$$g^\epsilon := dr^2 + \frac{1}{\epsilon^2}\sinh^2(\epsilon r)d\theta^2, \text{ and}$$
$$\overline{g}^\epsilon := dr^2 + \frac{1}{\epsilon^2}\sinh^2(\epsilon r)d\theta^2 + \cosh^2(\epsilon r)dt^2. \tag{3.1}$$

It is straightforward to show that this family varies smoothly with $\epsilon$. Furthermore, comparing with (2.3) and (2.4), we see that, for all $\epsilon > 0$, $g^\epsilon$ and $\overline{g}^\epsilon$ are related to $g$ and $\overline{g}$ by

$$g^\epsilon = \frac{1}{\epsilon^2}M_\epsilon^* g, \text{ and}$$
$$\overline{g}^\epsilon = \frac{1}{\epsilon^2}M_\epsilon^* \overline{g},$$

where $M_\epsilon$ here denotes uniform radial contraction by a factor of $\epsilon$ about the origin. It follows that, for all $\epsilon$, both $g^\epsilon$ and $\overline{g}^\epsilon$ are complete metrics of constant sectional curvature equal to $-\epsilon^2$.

We define a *euclidian minimal end* to be a surface $\Sigma$ in $\mathbb{R}^3$ which is a graph of some function $u$ over some annulus $A(R_0, \infty)$, which is minimal with respect to the metric $\overline{g}^0$ and which is symmetric by reflection in at least two vertical planes containing the $z$-axis. By the classical theory of minimal surfaces, the function $u$ then satisfies

$$u = a + c\text{Log}(r) + \text{O}\left(\frac{1}{r^{2+k}}\right). \tag{3.2}$$

for some real constants $a$ and $c$, which we call respectively the *constant term* and the *logarithmic parameter* of the euclidian minimal end*.

For $\epsilon > 0$, we define an *$\epsilon$-hyperbolic minimal end* to be a rotationally symmetric surface $\Sigma^\epsilon$ in $\mathbb{R}^3$ which is a graph of some function $u$ over some annulus $A(R_\epsilon, \infty)$ and which is minimal with respect to the metric $\overline{g}^\epsilon$. Given real constants $a$ and $c$, we set

$$u^\epsilon(R/4) = a + c\text{Log}(R/4),$$

and we set the flux $F$ equal to $2\pi c$. Upon rescaling and integrating (2.6), it then follows that, over the annulus $A(R/4, 2R^4)$

$$u^\epsilon = a + c\text{Log}(r) + \text{O}\left(r\left(\epsilon r + \frac{1}{r}\right)^3 r^{-k}\right). \tag{3.3}$$

---

\* In particular, planar ends will be considered as catenoidal ends with vanishing logarithmic parameters.

As before, we will refer to $a$ and $c$ respectively as the *constant term* and the *logarithmic parameter* of the $\epsilon$-hyperbolic end.

Now choose $(\epsilon, R)$ satisfying (2.1). Let $\Sigma$, $\Sigma^\epsilon$, $u$ and $u^\epsilon$ be as above and suppose that both surfaces have the same constant terms and logarithmic parameters. The $(\epsilon, R)$-*joined end* $\Sigma^{\epsilon, R}$ of $\Sigma$ is now defined to be the graph of the function

$$u^{\epsilon, R} := \chi_c u + (1 - \chi_c) u^\epsilon,$$

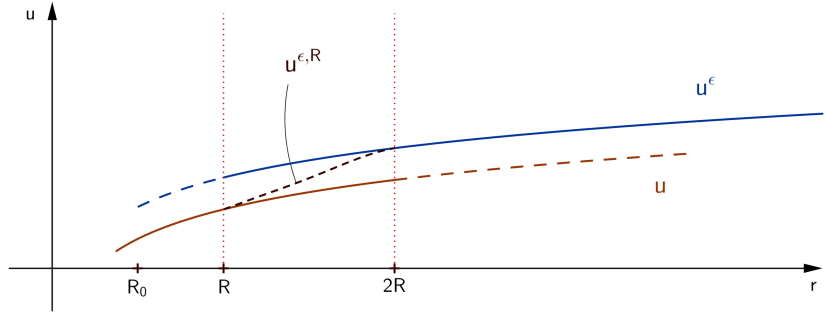

Figure 3.1.1 - The glueing operation - The graph of $u^{\epsilon, R}$ transitions from $u$ to $u^\epsilon$ over the annulus $A(R, 2R)$.

where $\chi_c$ is the cut-off function of the transition region $A(R, 2R)$ (c.f. Figure 3.1.1). Observe that

$$\begin{aligned}
\Sigma^{\epsilon, R} \cap (A(R_0, R) \times \mathbb{R}) &= \Sigma \cap (A(R_0, R) \times \mathbb{R}), \text{ and} \\
\Sigma^{\epsilon, R} \cap (A(2R, \infty) \times \mathbb{R}) &= \Sigma^\epsilon \cap (A(2R, \infty) \times \mathbb{R}),
\end{aligned} \tag{3.4}$$

whilst, over the annulus $A(R/4, 2R^4)$,

$$u^{\epsilon, R} = a + c \mathrm{Log}(r) + \mathrm{O}\left( r \left( \epsilon r + \frac{1}{r} \right)^3 \frac{1}{r^k} \right). \tag{3.5}$$

Finally, let $\Sigma$ be a complete, immersed surface in $\mathbb{R}^3$ such that $\Sigma \cap (B(R_0) \times \mathbb{R})$ is compact and $\Sigma \cap (A(R_0, \infty) \times \mathbb{R})$ is a union of finitely many euclidian minimal ends. For constants $(\epsilon, R)$ satisfying (2.1), we define the $(\epsilon, R)$-*joined surface* of $\Sigma$ to be the surface $\Sigma^{\epsilon, R}$ obtained by replacing each of its ends with the corresponding $(\epsilon, R)$-joined end (c.f. Figure 3.1.2).

**3.2 - Perturbation families.** Consider now a complete, immersed surface $\Sigma$ in $\mathbb{R}^3$ such that $\Sigma \cap (B(R_0) \times \mathbb{R})$ is compact and $\Sigma \cap (A(R_0, \infty) \times \mathbb{R})$ is a union of finitely many (not necessarily minimal) graphs over the annulus $A(R_0, \infty)$. Let $n$ be the number of ends of $\Sigma$ and let $\chi_0$, $\chi_0'$, $\chi_\epsilon'$ and $\chi_\epsilon$ be the cut-off functions of the transition regions $A(R_0, 2R_0)$, $A(2R_0, 4R_0)$, $A(1/2\epsilon, 1/\epsilon)$ and $A(1/\epsilon, 2/\epsilon)$ respectively. Define the vector field $X^\epsilon$ over $\Sigma$ by

$$X^\epsilon := \pm(\chi_\epsilon - \chi_0) e_z + (1 - (\chi_\epsilon - \chi_0)) N^\epsilon, \tag{3.6}$$

Minimal surfaces in hyperbolic space.

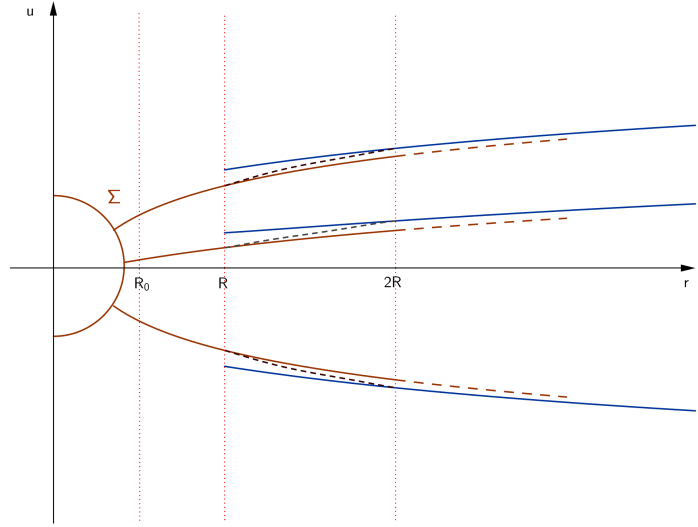

Figure 3.1.2 - The joined surface - Each of the ends of $\Sigma$ is replaced with its corresponding $(\epsilon,R)$-joined end. The non-trivial part of the topology of $\Sigma$, which is contained in the cylinder $B(R) \times \mathbb{R}$, is represented here schematically by a semi-circle.

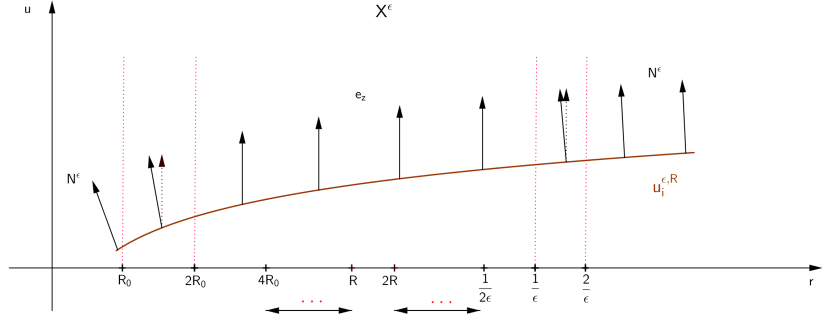

Figure 3.2.3 - The modified normal vector field - The modified normal vector field $X^\epsilon$ coincides with the unit hyperbolic normal over the complement of the annulus $A(R_0, 2/\epsilon)$, and coincides with the third basis vector over the annulus $A(2R_0, 1/\epsilon)$.

where $N^\epsilon$ here denotes the unit normal vector field over $\Sigma$ with respect to the metric $\overline{g}^\epsilon$ and the sign of $e_z$ is chosen over each end so that it points in the same direction as $N^\epsilon$ (c.f. Figure 3.2.3). Define $\mathcal{E}^\epsilon : \mathbb{R}^n \times \mathbb{R}^n \times C_0^\infty(\Sigma) \to C^\infty(\Sigma, \mathbb{R}^3)$ by

$$\mathcal{E}^\epsilon_{a,b,f}(x) := x + f(x)X^\epsilon(x) + \sum_{i=1}^{n} \mathbb{I}_i(x)\big(a_i\big(1 - \chi_0'(x)\big) + b_i\big(1 - \chi_\epsilon'(x)\big)\big)e_z, \qquad (3.7)$$

where, for all $1 \leq i \leq n$, $\mathbb{I}_i$ here denotes the indicator function of the $i$'th component of $\Sigma \cap (A(R_0, \infty) \times \mathbb{R})$. Observe that, for all $a$ and $b$ and for all sufficiently small $f$, $\mathcal{E}^\epsilon_{a,b,f}$ defines an immersion of $\Sigma$ into $\mathbb{R}^3$. Heuristically, $\mathcal{E}_{a,b,f}^\epsilon$ decomposes into a perturbation of $\Sigma$ in $f$ times the modified normal direction, followed by a vertical perturbation of each end by a height of $a_i$ outside the cylinder $B(4R_0) \times \mathbb{R}$, and then a further vertical perturbation of each end by a height of $b_i$ outside the cylinder $B(2/\epsilon) \times \mathbb{R}$ (c.f. Figures 3.2.4, 3.2.5 and 3.2.6). Although the second vertical perturbations may appear superfluous, they serve to control singular terms arising from the Green's operators of the $\epsilon$-hyperbolic ends, as we will see in Section 4.2, below.

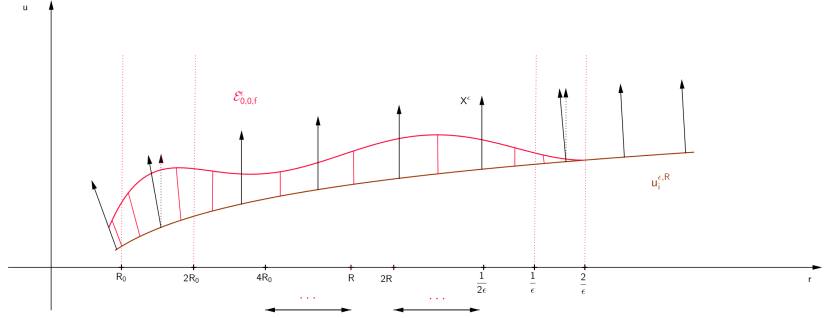

Figure 3.2.4 - A microscopic perturbation - The surface $\Sigma$ is perturbed in the direction of the modified normal vector field. Such perturbations are referred to as microscopic perturbations since, the function $f$ decaying at infinity, they may be thought of as being supported over the annulus $A(2\epsilon R^4)$, which becomes negligable as $\Lambda$ tends to infinity.

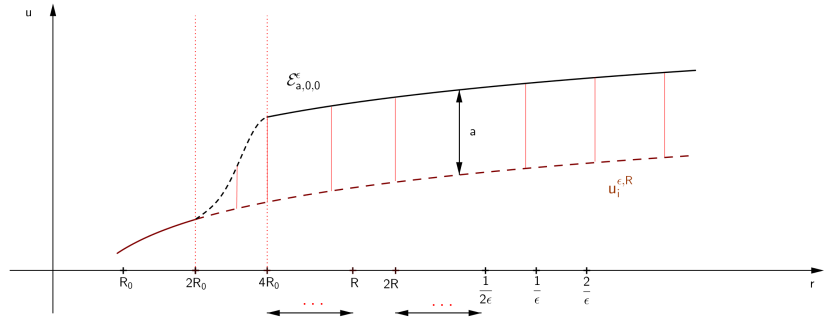

Figure 3.2.5 - The first macroscopic perturbation - The end of the surface $\Sigma$ is shifted vertically over the annulus $A(4R_0, \infty)$. Such perturbations are referred to as macroscopic perturbations, as they remain significant for all large values of $\Lambda$.

This family generates perturbation operators as follows. First, define $\mathcal{H}^\epsilon : \mathbb{R}^n \times \mathbb{R}^n \times C_0^\infty(\Sigma) \to C_0^\infty(\Sigma)$ such that $\mathcal{H}_{a,b,f}^\epsilon$ is the mean curvature function of the immersion $\mathcal{E}_{a,b,f}^\epsilon$ with respect to the metric $\overline{g}^\epsilon$. In the present context, we define the *modified Jacobi operator* of $\Sigma$ with respect to the metric $\overline{g}^\epsilon$ by

$$\hat{J}_\Sigma^\epsilon f := \frac{1}{\psi^\epsilon} \frac{\partial}{\partial t} \mathcal{H}_{0,0,tf}^\epsilon(x) \bigg|_{t=0}, \tag{3.8}$$

Minimal surfaces in hyperbolic space.

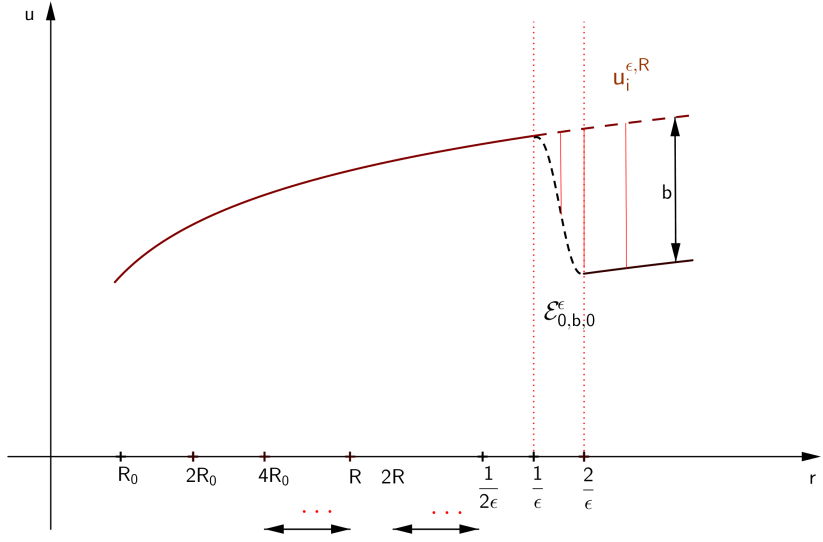

Figure 3.2.6 - The second macroscopic perturbation - The end of the surface $\Sigma$ is shifted vertically over the annulus $A(2/\epsilon, \infty)$. As before, such perturbations are also referred to as macroscopic perturbations.

where
$$\psi^\epsilon := \overline{g}^\epsilon\big(X^\epsilon, N^\epsilon\big).$$
In addition, we define the operators $Y_\Sigma^\epsilon, Z_\Sigma^\epsilon : \mathbb{R}^n \to C_0^\infty(\Sigma)$ by

$$\begin{aligned}
(Y_\Sigma^\epsilon v)(x) &:= \frac{1}{\psi^\epsilon}\left.\frac{\partial}{\partial t}\mathcal{H}_{tv,0,0}^\epsilon(x)\right|_{t=0}, \text{ and}\\
(Z_\Sigma^\epsilon w)(x) &:= \frac{1}{\psi^\epsilon}\left.\frac{\partial}{\partial t}\mathcal{H}_{0,tw,0}^\epsilon(x)\right|_{t=0}.
\end{aligned} \tag{3.9}$$

Observe that the above construction remains valid in the case where $\epsilon = 0$ where, by convention, the cut-off functions $\chi_\epsilon$ and $\chi'_\epsilon$ are taken to be identically equal to 1.

In the case where each end of $\Sigma$ satisfies (3.2) for some constants $a$ and $c$, we also consider the following, slightly more subtle, family of perturbations. For each $1 \leq i \leq n$, let $a_{0,i}$ and $c_{0,i}$ be respectively the constant term and the logarithmic parameter of the $i$'th end of $\Sigma$. Let $(\Sigma_c)_{c \in \mathbb{R}^n}$ be a smoothly varying family of immersed surfaces in $\mathbb{R}^3$ such that $\Sigma_{c_0} = \Sigma$ and such that, for all $c$, and for all $1 \leq i \leq n$, the $i$'th component of $\Sigma_c \cap (A(R_0, \infty) \times \mathbb{R})$ satisfies (3.2) with constant term $a_{0,i}$ and logarithmic parameter $c_i$. Now let $\mathcal{F} : \mathbb{R}^n \times \Sigma \to \mathbb{R}^3$ be a smooth function such that

(1) for all $c$, $\mathcal{F}_c := \mathcal{F}(c, \cdot)$ parametrises $\Sigma_c$; and

(2) for all $c$, and for all $x \in \Sigma \cap (A(R_0, +\infty) \times \mathbb{R})$,

$$(\pi \circ \mathcal{F}_c)(x) = \pi(x),$$

Minimal surfaces in hyperbolic space.

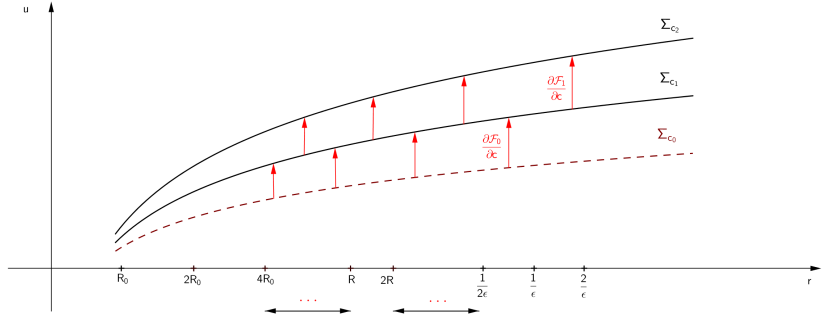

Figure 3.2.7 - The third macroscopic perturbation - The logarithmic parameter of
the end of the surface $\Sigma$ is continuously varied.

where $\pi : \mathbb{R}^3 \rightarrow \mathbb{R}^2$ is the canonical projection onto the first two coordinates (c.f. Figure 3.2.7).

We may suppose that, for all $c$, $\mathcal{F}_c$ defines an immersion of $\Sigma$ into $\mathbb{R}^3$. Depending on the context, we will also impose further conditions on the family $(\Sigma_c)_{c \in \mathbb{R}^n}$ used to define $\mathcal{F}$. Indeed, if the surface $\Sigma$ is known to possess certain symmetries, then we will assume that, for all $c$, $\Sigma_c$ also possesses the same symmetries and if $\Sigma$ is a Costa-Hoffman-Meeks surface, as in Section 4.1, below, then we will assume that, for all $c$, every end of $\Sigma_c$ is a euclidian minimal end.

A similar construction also exists in the case where each of the ends of $\Sigma$ are $(\epsilon, R)$-joined ends. Here the family $(\Sigma_c)$ is chosen such that, for all $c$, every end of $\Sigma_c$ is also an $(\epsilon, R)$-joined end. The rest of the construction of $\mathcal{F}$ continues as before.

In each case, the family $\mathcal{F}$ generates an additional perturbation operator as follows. Indeed, for $\epsilon \geq 0$, we define $\mathcal{H}^\epsilon : \mathbb{R}^n \rightarrow C^\infty(\Sigma)$ such that $\mathcal{H}_c^\epsilon$ is the mean curvature of the immersion $\mathcal{F}_c$ with respect to the metric $\overline{g}^\epsilon$. We define the operator $X_\Sigma^\epsilon : \mathbb{R}^n \rightarrow C_0^\infty(\Sigma)$ by

$$(X_\Sigma^\epsilon u)(x) := \frac{1}{\psi^\epsilon} \left. \frac{\partial}{\partial t} \mathcal{H}_{c_0 + tu}^{\epsilon, R}(x) \right|_{t=0}, \tag{3.10}$$

where $\psi^\epsilon$ is the function defined above.

Finally, the above two families are combined as follows. First, in order to emphasize the dependence of $\mathcal{E}^\epsilon$ on the initial surface $\Sigma$, we denote

$$\mathcal{E}^\epsilon := \mathcal{E}^\epsilon[\Sigma].$$

Next, for all $c \in \mathbb{R}^n$, we identify functions over $\Sigma_c$ with functions over $\Sigma$ by composition with $\mathcal{F}_c$. We now define $\tilde{\mathcal{E}}^\epsilon : \mathbb{R}^n \times \mathbb{R}^n \times \mathbb{R}^n \times C_0^\infty(\Sigma) \rightarrow C^\infty(\Sigma, \mathbb{R}^3)$ by

$$\tilde{\mathcal{E}}_{c,a,b,f}^{\epsilon, R} := \mathcal{E}_{a,b,f}^\epsilon[\Sigma_c] \circ \mathcal{F}_c.$$

Theorem A will be proven by constructing such a family around the $(\epsilon, R)$-joined surface of a Costa-Hoffman-Meeks surface and by showing that, for sufficiently large $\Lambda$, there is a point $(c, a, b, f)$ close to $(c_0, 0, 0, 0)$ such that the immersion $\tilde{\mathcal{E}}_{c,a,b,f}^{\epsilon, R}$ is minimal with respect to $\overline{g}^\epsilon$.

**3.3 - The perturbation operators.** We conclude this chapter by reviewing certain functional properties of the operators constructed above. We continue to use the notation and terminology of the previous sections.

**Lemma 3.3.1**

*If $\Sigma$ is a euclidian minimal end, an $\epsilon$-hyperbolic minimal end, or an $(\epsilon, R)$-joined end, and if $\hat{J}_\Sigma^\epsilon$ is its modified Jacobi operator with respect to the metric $\overline{g}^\epsilon$ then, over the annulus $A(R/4, 4R^4)$,*

$$\hat{J}_\Sigma^\epsilon = \Delta^\epsilon f - \frac{c^2 x^i x^j}{r^4} f_{ij} + 2\left(\epsilon^2 + \frac{c^2}{r^4}\right) x^i f_i + \mathcal{R}, \tag{3.11}$$

*where $\Delta^\epsilon$ here denotes the Laplace-Beltrami operator of the metric $g^\epsilon$, $r$ denotes the radial distance from the origin, and the remainder $\mathcal{R}$ is of the form*

$$\mathcal{R}f = a^{ij} f_{ij} + b^i f_i, \tag{3.12}$$

*where the coefficients $a^{ij}$ and $b^i$ satisfy*

$$a^{ij} = O\left(\left(\epsilon r + \frac{1}{r}\right)^4 \frac{1}{r^k}\right), \text{ and}$$
$$b^i = O\left(\left(\epsilon r + \frac{1}{r}\right)^4 \frac{1}{r^{k+1}}\right). \tag{3.13}$$

**Proof:** Indeed, in each case, by (3.2), (3.3) and (3.5), over the annulus $A(R/4, 2R^4)$,

$$u = a + c\text{Log}(r) + \text{O}\left(r\left(\epsilon r + \frac{1}{r}\right)^3 \frac{1}{r^k}\right),$$

and the result follows upon rescaling (2.16). $\square$

**Lemma 3.3.2**

*If $\Sigma$ is a euclidian minimal end and if, for all $\epsilon$, $\hat{J}_\Sigma^\epsilon$ is its modified Jacobi operator with respect to the metric $\overline{g}^\epsilon$ then, over the annulus $A(R_0, 2R^4)$,*

$$\left(\hat{J}_\Sigma^\epsilon - \hat{J}_\Sigma^0\right) = a^{ij} f_{ij} + b^i f_i, \tag{3.14}$$

*where the coefficients $a^{ij}$ and $b^i$ satisfy*

$$a^{ij} = O\left(\left(\epsilon^2 r^2 + \frac{\epsilon}{r^2} + \frac{1}{r^4}\right)\frac{1}{r^k}\right), \text{ and}$$
$$b^i = O\left(\left(\epsilon^2 r^2 + \frac{\epsilon}{r^2} + \frac{1}{r^4}\right)\frac{1}{r^{k+1}}\right), \tag{3.15}$$

**Proof:** Indeed, for all $\epsilon \geq 0$, let $\Delta^\epsilon$ denote the Laplace-Beltrami operator of the metric $g^\epsilon$. We readily verify that

$$\Delta^\epsilon = \frac{(\epsilon r)^2}{\sinh^2(\epsilon r)} f_{\theta\theta} + \epsilon\coth(\epsilon r) f_r + f_{rr}, \tag{3.16}$$

where $f_r$ and $f_\theta$ here denote respectively the derivatives of $f$ in the unit radial and unit angular directions with respect to the euclidean metric over $\mathbb{R}^2$. It follows that, over $A(R_0, 2R^4)$,

$$\Delta^\epsilon f = \Delta^0 f + a^{ij} f_{ij} + b^i f_i,$$

where

$$a^{ij} = O\left(\frac{\epsilon^2 r^2}{r^k}\right), \text{ and}$$

$$b^i = O\left(\frac{\epsilon^2 r}{r^k}\right).$$

The result now follows by (2.1) and (3.11). $\square$

Finally, it is straightforward to show that the supports of $X_\Sigma^\epsilon$, $Y_\Sigma^\epsilon$ and $Z_\Sigma^\epsilon$ satisfy

$$\mathrm{Supp}(X_\Sigma^\epsilon) \subseteq \Sigma \cap (A(2R_0, 4R_0) \times \mathbb{R}),$$
$$\mathrm{Supp}(Y_\Sigma^\epsilon) \subseteq \Sigma \cap (B(2R) \times \mathbb{R}), \text{ and}$$
$$\mathrm{Supp}(Z_\Sigma^\epsilon) \subseteq \Sigma \cap (A(1/2\epsilon, 1/\epsilon) \times \mathbb{R}).$$

Furthermore, for all $v, w \in \mathbb{R}^n$, and for all $x \in \Sigma$,

$$\begin{aligned}
(Y_\Sigma^\epsilon v)(x) &= -\sum_{i=1}^n \mathbb{I}_i(x) v_i (\hat{J}_\Sigma^\epsilon \chi_0')(p), \text{ and} \\
(Z_\Sigma^\epsilon w)(x) &= -\sum_{i=1}^n \mathbb{I}_i(x) w_i (\hat{J}_\Sigma^\epsilon \chi_\epsilon')(p).
\end{aligned} \tag{3.17}$$

It thus only remains to estimate $X_\Sigma^\epsilon$.

**Lemma 3.3.3**

If $\Sigma$ is an $(\epsilon, R)$-joined end, then, over $A(2R_0, R)$,

$$X_\Sigma^\epsilon u = O\left(\frac{1}{r^{6+k}}\right), \tag{3.18}$$

and over $A(R, 2R)$,

$$X_\Sigma^\epsilon u = O\left(\frac{1}{r^{4+k}}\right). \tag{3.19}$$

**Proof:** Let $\Sigma_c$ be a smooth family of ends as in Section 3.2. For all $t$, let $f_t : A(2R_0, \infty) \to \mathbb{R}$ denote the profile of $\Sigma_{c_0+t}$ and define

$$\phi := \left.\frac{\partial f_t}{\partial t}\right|_{t=0}.$$

Over $A(2R_0, 2R)$, we have

$$X_\Sigma^\epsilon u = \hat{J}_\Sigma^\epsilon \phi.$$

Furthermore, rescaling (2.7) yields

$$\phi = \mathrm{Log}(r) + \mathrm{O}\left(\frac{1}{r^{2+k}}\right), \tag{3.20}.$$

By hypothesis, $\Sigma_c \cap (A(2R_0, R) \times \mathbb{R})$ is minimal with respect to $\bar{g}^0$ for all $c$ so that

$$\hat{J}_\Sigma^0 \phi = 0.$$

The first estimate now follows from (3.14) and (2.1). Finally, by (3.16),

$$\Delta^\epsilon \mathrm{Log}(r) = \mathrm{O}\left(\frac{\epsilon^2}{r^k}\right).$$

The second estimate now follows by (3.11) and (2.1), and this completes the proof. $\square$

# 4 - Constructing the Green's operator.

**4.1 - Pertubation theory.** Let $\Sigma$ be the Costa-Hoffman-Meeks surface of genus $g$ for some positive integer $g$ and let $\hat{\Sigma}$ be its $(\epsilon, R)$-joined surface, as constructed in Section 3.1, for some $(\epsilon, R)$ satisfying (2.1). Observe that there is an (almost) natural diffeomorphism from $\Sigma$ into $\hat{\Sigma}$ which maps points in each of the ends of $\Sigma$ directly upwards or downwards. Consequently, functions over $\Sigma$ may equally well be considered as functions over $\hat{\Sigma}$. Similarly, functions over $\hat{\Sigma} \cap (A(R_0, \infty) \times \mathbb{R})$ may be considered as functions over three copies of $A(R_0, \infty)$, and so on. In defining norms over spaces of functions, we will pass between these different perspectives without comment.

We now study operators and function spaces defined over these surfaces. First, define the *scale-free total derivative operator* by

$$D_{\mathrm{SF}} := rD, \tag{4.1}$$

where $D$ here denotes the total derivative operator of $\mathbb{R}^2$ and $r$ denotes the radial distance from the origin. For $\alpha \in [0, 1]$, define the *scale free Hölder operator* by

$$\delta_{\mathrm{SF}}^\alpha f(r) := r^\alpha \left[f|_{A(r/2, 2r)}\right]_\alpha. \tag{4.2}$$

For all non-negative integer $m$, for all $\alpha \in [0,1]$, for all real $\delta$, define the *scale-free weighted Hölder norm* of weight $\delta$ over $A(R_0, \infty)$ such that, for every $m$-times differentiable function $f : A(R_0, \infty) \to \mathbb{R}$,

$$\|f\|_{C^{m,\alpha}_{\delta,\mathrm{SF}}(A(R_0,\infty))} := \sum_{k=0}^{\infty} \left\| r^\delta D^k_{\mathrm{SF}} f \right\|_{C^0(A(R_0,\infty))} + \left\| r^\delta \delta^\alpha_{\mathrm{SF}} D^n_{\mathrm{SF}} f \right\|_{C^0(]2R_0,\infty[)}, \qquad (4.3)$$

and define the scale-free weighted Hölder norm of weight $\delta$ over $\Sigma$ such that, for every $m$-times differentiable function $f : \Sigma \to \mathbb{R}$,

$$\|f\|_{C^{m,\alpha}_{\delta,\mathrm{SF}}(\Sigma)} := \left\| f|_{\Sigma \cap (B(2R_0) \times \mathbb{R})} \right\|_{C^{m,\alpha}} + \left\| f|_{\Sigma \cap (A(R_0,\infty) \times \mathbb{R})} \right\|_{C^{m,\alpha}_{\delta,\mathrm{SF}}(A(R_0,\infty))}. \qquad (4.4)$$

For all such $m$, $\alpha$ and $\delta$, let $C^{m,\alpha}_{\delta,\mathrm{SF},g}(\Sigma)$ denote the space of $m$-times differentiable functions $f$ over $\Sigma$ which satisfy $\|f\|_{C^{m,\alpha}_{\delta,\mathrm{SF}}(\Sigma)} < \infty$ and which also satisfy $f \circ \sigma = f$ for every horizontal symmetry $\sigma$ of $\Sigma$.

Consider now the operators $\hat{J}^0_\Sigma$, $X^0_\Sigma$ and $Y^0_\Sigma$ defined over the Costa-Hoffman-Meeks surface $\Sigma$ as in Section 3.2. Observe that, in this case, the operator $X^0_\Sigma$ is defined using a family $(\Sigma_c)_{c \in \mathbb{R}^3}$ chosen such that, for all $c$, $\Sigma_c$ is a complete immersed surface, symmetric under all the horizontal symmetries of the Costa-Hoffman-Meeks surface, all of whose ends are euclidian minimal ends. From the perturbation theory of Costa-Hoffman-Meeks surfaces (c.f. [6], [13], [14] and [17]), we have

**Theorem 4.1.1**

*For all $\alpha \in ]0,1[$, for all $\delta \in ]1,2[$, and for all $R_0$ sufficiently large, $(X^0_\Sigma, Y^0_\Sigma, \hat{J}^0_\Sigma)$ defines a surjective Fredholm map of Fredholm index 3 from $\mathbb{R}^3 \oplus \mathbb{R}^3 \oplus C^{2,\alpha}_{\delta,\mathrm{SF},g}(\Sigma)$ into $C^{0,\alpha}_{\delta+2,\mathrm{SF},g}(\Sigma)$.*

For $\epsilon > 0$, consider now the operators $\hat{J}^\epsilon_{\hat{\Sigma}}$, $X^\epsilon_{\hat{\Sigma}}$, $Y^\epsilon_{\hat{\Sigma}}$ and $Z^\epsilon_{\hat{\Sigma}}$ defined over the $(\epsilon, R)$-joined surface $\hat{\Sigma}$ as in Section 3.2. Observe that, this time, the operator $X^\epsilon_{\hat{\Sigma}}$ is defined using a family $(\hat{\Sigma}_c)_{c \in \mathbb{R}^3}$ chosen such that, for all $c$, $\hat{\Sigma}_c$ is a complete immersed surface, symmetric under all the horizontal symmetries of the Costa-Hoffman-Meeks surface, all of whose ends are $(\epsilon, R)$-joined ends. For technical reasons, it will be useful at this stage to alter slightly the modified Jacobi operator of the Costa-Hoffman-Meeks surface, and we define

$$\hat{J}^{\epsilon,1}_\Sigma := \chi_{2R^4} \hat{J}^\epsilon_\Sigma + (1 - \chi_{2R^4}) \hat{J}^0_\Sigma,$$

where $\chi_{2R^4}$ here denotes the cut-off function of the transition region $A(2R^4, 4R^4)$. In particular, this operator coincides with $\hat{J}^\epsilon_{\hat{\Sigma}}$ over the region $\Sigma \cap (B(2R^4) \times \mathbb{R})$, which is where it will be explicitely used.

**Theorem 4.1.2**

*For all $\alpha \in ]0,1[$, for all $\delta \in ]1,2[$, for all $R_0 > 0$ sufficiently large, and for all sufficiently large $\Lambda$, the triplet $(X^\epsilon_{\hat{\Sigma}}, Y^\epsilon_{\hat{\Sigma}}, \hat{J}^{\epsilon,1}_\Sigma)$ defines a surjective Fredholm map from $\mathbb{R}^3 \oplus \mathbb{R}^3 \oplus C^{2,\alpha}_{\delta,\mathrm{SF},g}(\Sigma)$ into $C^{0,\alpha}_{2+\delta,\mathrm{SF},g}(\Sigma)$ of Fredholm index 3. Furthermore, the right inverse $(U, V, \Phi)$ can be chosen in such a manner that its norm is uniformly bounded, independent of $\Lambda$.*

**Remark:** In the sequel, $R_0$ will be chosen large enough for Lemma 4.1.2 to hold for all large values of $\Lambda$. It will then be fixed once and for all, and $\Lambda$ will be made to tend to $+\infty$.

**Proof:** Since both $Y_{\hat{\Sigma}}^{\epsilon}$ and $Y_{\Sigma}^{0}$ are supported over $\Sigma \cap (B(2R_0) \times \mathbb{R})$, and since $\overline{g}^{\epsilon}$ converges to $\overline{g}^{0}$ as $\Lambda$ tends to infinity, it follows that the difference $(Y_{\hat{\Sigma}}^{\epsilon} - Y_{\Sigma}^{0})$, considered as a map from $\mathbb{R}^3$ into $C_{2+\delta, \text{SF}, g}^{0,\alpha}(\Sigma)$, tends to 0 as $\Lambda$ tends to infinity. Likewise, by (3.14) and (3.15), the difference $(\hat{J}_{\Sigma}^{\epsilon,1} - \hat{J}_{\Sigma}^{0})$, considered as a map from $C_{\delta,\text{SF},g}^{2,\alpha}(\Sigma)$ into $C_{\delta+2,\text{SF},g}^{0,\alpha}(\Sigma)$, also tends to zero as $\Lambda$ tends to infinity. It remains to estimate the difference $(X_{\hat{\Sigma}}^{\epsilon} - X_{\Sigma}^{0})$. Observe first that this operator is supported over $\Sigma \cap (B(2R) \times \mathbb{R})$. Next, as above, for all $R_0$,

$$\left\| (X_{\hat{\Sigma}}^{\epsilon} - X_{\Sigma}^{0})|_{\Sigma \cap (B(2R_0) \times \mathbb{R})} \right\|_{C^{0,\alpha}} \to 0$$

as $\Lambda$ tends to infinity. Finally, observe that $X_{\Sigma}^{0}$ vanishes over $\Sigma \cap (A(R_0, 2R) \times \mathbb{R})$ so that, by (3.18) and (3.19),

$$\left\| (X_{\hat{\Sigma}}^{\epsilon} - X_{\Sigma}^{0})|_{\Sigma \cap (A(R_0,2R) \times \mathbb{R})} \right\|_{C_{\delta+2,\text{SF}}^{0,\alpha}} = \left\| X_{\hat{\Sigma}}^{\epsilon}|_{\Sigma \cap (A(R_0,2R) \times \mathbb{R})} \right\|_{C_{\delta+2,\text{SF}}^{0,\alpha}} \to 0$$

as $\Lambda$ tends to infinity. The result now follows by stability of the surjectivity property for Fredholm maps. $\square$

For all non-negative, integer $m$, for all $\alpha \in [0,1]$, for all $\gamma \in \mathbb{R}$ and for all $\epsilon > 0$, define the following *weighted Sobolev and Hölder norms* for functions over $\mathbb{R}^2$,

$$\begin{aligned} \|f\|_{H_{\gamma,\epsilon}^m(\mathbb{R}^2)} &:= \|f(\cdot/\epsilon)\|_{H_{\gamma}^m(\mathbb{H}^2)}. \text{ and} \\ \|f\|_{C_{\gamma,\epsilon}^{m,\alpha}(\mathbb{R}^2)} &:= \|f(\cdot/\epsilon)\|_{C_{\gamma}^{m,\alpha}(\mathbb{H}^2)}, \end{aligned} \tag{4.5}$$

and define the hybrid norm by

$$\|f\|_{m,\alpha,\gamma,\epsilon} := \|f\|_{C_{\gamma,\epsilon}^{m,\alpha}(\mathbb{R}^2)} + \frac{1}{\epsilon R} \|f\|_{H_{\gamma,\epsilon}^m(\mathbb{R}^2)}. \tag{4.6}$$

For all such $m$, $\alpha$, $\gamma$ and $\epsilon$, let $\mathcal{L}_{\gamma,\epsilon,g}^{m,\alpha}(\mathbb{R}^2)$ denote the space of $m$-times differentiable functions which satisfy $\|f\|_{m,\alpha,\gamma,\epsilon}$, and which also satisfy $f \circ \sigma = f$ for every horizontal symmetry $\sigma$ of the Costa-Hoffman-Meeks surface.

Let $\hat{\Sigma}^{\epsilon}$ be an $\epsilon$-hyperbolic end over the annulus $A(R/4, +\infty)$, and let $\hat{J}_{\hat{\Sigma}_{\epsilon}}$ be its modified Jacobi operator, as defined in Section 2.2. Upon rescaling, Theorem 2.3.3 immediately yields

### Lemma 4.1.3

*For sufficiently small $\alpha \in ]0,1[$, for sufficiently small $\gamma$, and for sufficiently large $\Lambda$, the operator $\hat{J}_{\hat{\Sigma}_{\epsilon}}^{\epsilon}$ defines a linear isomorphism from $\mathcal{L}_{\gamma,\epsilon,g}^{2,\alpha}(\mathbb{R}^2)$ into $\mathcal{L}_{\gamma,\epsilon,g}^{0,\alpha}(\mathbb{R}^2)$. Furthermore, denoting its inverse by $\Psi$, we may suppose that the operator norm of $\epsilon^2 \Psi$ is uniformly bounded independent of $\Lambda$.*

**4.2 - The ping-pong argument.** For notational convenience, we henceforth work as if $\Sigma$ and $\hat\Sigma$ each had only one end. Consider now the following seminorms for functions over $\hat\Sigma$.

$$\|f\|_{m,E} := \|f|_{B(0,4R)}\|_{C^{m,\alpha}_{(2-m)+\delta,\mathrm{SF}}(\Sigma)},$$

$$\|f\|_{m,F,\mathrm{H\ddot ol}} := \|f|_{A(R,\infty)}\|_{C^{m,\alpha}_{\gamma,\epsilon}(\mathbb{R}^2)},$$

$$\|f\|_{m,F,\mathrm{Sob}} := \|f|_{A(R,\infty)}\|_{H^m_{\gamma,\epsilon}(\mathbb{R}^2)},\ \text{and} \tag{4.7}$$

$$\|f\|_{m,F} := \|f\|_{m,F,\mathrm{H\ddot ol}} + \frac{1}{\epsilon R}\|f\|_{m,F,\mathrm{Sob}}.$$

Let $\mathcal{E}$ denote the space of continuous functions with support in $\hat\Sigma \cap (B(4R)\times\mathbb{R})$ which have finite $\|\cdot\|_{0,E}$-norm and which are invariant under every horizontal symmetry of the Costa-Hoffman-Meeks surface. Likewise, let $\mathcal{F}$ denote the space of continuous functions with support in $\hat\Sigma \cap (A(R,\infty)\times\mathbb{R})$ which have finite $\|\cdot\|_{0,F}$-norm and which are also invariant under these symmetries. Define the operator $A : \mathcal{E} \to \mathcal{F}$ by

$$Ae := \hat{J}^\epsilon_{\hat\Sigma}\chi_u\Phi e + X^\epsilon_{\hat\Sigma}Ue + Y^\epsilon_{\hat\Sigma}Ve - e, \tag{4.8}$$

where $\chi_u$ is the cut-off function of the transition region $A(R^4,2R^4)$ and $(U,V,\Phi)$ is defined as in Lemma 4.1.2. This operator measures the extent to which $(U,V,\chi_u\Phi)$ fails to be a Green's operator of $(X^\epsilon_{\hat\Sigma},Y^\epsilon_{\hat\Sigma},\hat{J}^\epsilon_{\hat\Sigma})$ for functions in $\mathcal{E}$. In particular, since $\hat{J}^\epsilon_{\hat\Sigma}$ coincides with $\hat{J}^{\epsilon,1}_{\hat\Sigma}$ over $B(0,R)$, $Ae$ is supported in the interior of $A(R,\infty)$ making it indeed an element of $\mathcal{F}$.

Likewise, define the operators $B : \mathcal{F} \to \mathcal{E}$ and $W : \mathcal{F} \to \mathbb{R}^3$ by,

$$Bf := \hat{J}^\epsilon_{\hat\Sigma}(1-\chi_l)(\Psi f - \chi'_\epsilon(Wf)) - Z^\epsilon_{\hat\Sigma}Wf - f,\ \text{and}$$
$$Wf := (\Psi f)(0). \tag{4.9}$$

where $\chi_l$ and $\chi'_\epsilon$ are the cut-off functions of the transition regions $A(R/4,R/2)$ and $A(1/2\epsilon,1/\epsilon)$ and $\Psi$ is defined as in Lemma 4.1.3. This operator measures the extent to which $(-W,(1-\chi_l)(\Psi-\chi'_\epsilon W))$ fails to be a Green's operator of $(Z^\epsilon_{\hat\Sigma},\hat{J}^\epsilon_{\hat\Sigma})$ for functions in $\mathcal{F}$. In particular, by (3.17) together with the fact that $\hat{J}^\epsilon_{\hat\Sigma}$ coincides with $\hat{J}^\epsilon_\Sigma$ over $A(2R,\infty)$, $Bf$ is supported in $B(4R)$, and is thus indeed an element of $\mathcal{E}$.

In Section 6 of [18], we prove

**Theorem 4.2.1**

*For all $\delta \in ]1,2[$, for sufficiently small $\alpha$, for all sufficiently small $\gamma$, for all $e \in \mathcal{E}$ and for all $f \in \mathcal{F}$,*

$$\|Ae\|_{0,F} \lesssim \frac{1}{(\epsilon R)^{2\alpha}}\frac{1}{R^{6+\delta}}\|e\|_{0,E}.\ \text{and}$$

$$\|Bf\|_{0,E} \lesssim \frac{R^2}{(\epsilon R)}\|f\|_{0,F}. \tag{4.10}$$

**Remark:** The proof of Theorem 4.2.1 consists of straightforward, though highly technical, estimates of every one of the terms involved. The context of [18] is slightly different to the one studied here. However, this merely simplifies our calculations in the present case, as the operators $(\hat{J}_{\hat{\Sigma}}^{\epsilon} - \hat{J}_{\Sigma}^{\epsilon '})$ and $(\hat{J}_{\hat{\Sigma}}^{\epsilon} - \hat{J}_{\Sigma^{\epsilon}}^{\epsilon})$ here contain fewer terms and all other terms are unaltered.

Bearing in mind (2.1), it follows from Theorem 4.2.1 that, for $\delta \in ]1, 2[$, for sufficiently small $\alpha$ and for sufficiently small $\gamma$, the operator norms of the compositions $AB$ and $BA$ satisfy

$$\|AB\|, \|BA\| \lesssim \frac{1}{(\epsilon R)^{2\alpha}} \frac{1}{\epsilon R^{5+\delta}} \lesssim \frac{1}{\Lambda} < 1.$$

We therefore define $Q_E : \mathcal{E} \to \mathcal{E}$ and $Q_F : \mathcal{F} \to \mathcal{F}$ by

$$Q_E := \sum_{m=0}^{\infty} (BA)^m, \text{ and}$$
$$Q_F := \sum_{m=0}^{\infty} (AB)^m, \tag{4.11}$$

and we define

$$\hat{U}f := UQ_E\chi f - UBQ_F(1-\chi)f,$$
$$\hat{V}f := VQ_E\chi f - VBQ_F(1-\chi)f,$$
$$\hat{W}f := WAQ_F\chi f - WQ_F(1-\chi)f, \text{ and} \tag{4.12}$$
$$\hat{P}f := P_F\chi f + P_G(1-\chi)f,$$

where $\chi$ is the cut-off function of the transition region $A(2R, 4R)$ and

$$P_F e := \chi_u \Phi Q_E e - (1 - \chi_l)\big(\Psi A Q_E e - \chi_\epsilon'(W A Q_E e)\big), \text{ and}$$
$$P_G f := -\chi_u \Phi B Q_F f + (1 - \chi_l)\big(\Psi Q_F f - \chi_\epsilon'(W Q_F f)\big). \tag{4.13}$$

In Section 6.4 of [18], we show that, for all suitable $f$,

$$\hat{J}_{\hat{\Sigma}}^{\epsilon}\hat{P}f + X_{\hat{\Sigma}}^{\epsilon}\hat{U}f + Y_{\hat{\Sigma}}^{\epsilon}\hat{V}f + Z_{\hat{\Sigma}}^{\epsilon}\hat{W}f = f, \tag{4.14}$$

so that $(\hat{U}, \hat{V}, \hat{W}, \hat{P})$ is a Green's operator for $(X_{\hat{\Sigma}}^{\epsilon}, Y_{\hat{\Sigma}}^{\epsilon}, Z_{\hat{\Sigma}}^{\epsilon}, \hat{J}_{\hat{\Sigma}}^{\epsilon})$. Furthermore, we obtain

**Theorem 4.2.2**

*For all $\delta \in ]1, 2[$, for all sufficiently small $\alpha \in ]0, 1[$, for all sufficiently small $\gamma$, and for all $f$,*

$$\|\hat{U}f\| \lesssim \|f\|_{0,E} + \frac{R^2}{(\epsilon R)^{1+\alpha}}\|f\|_{0,F},$$

$$\|\hat{V}f\| \lesssim \|f\|_{0,E} + \frac{R^2}{(\epsilon R)^{1+\alpha}}\|f\|_{0,F},$$

$$\|\hat{W}f\| \lesssim \|f\|_{0,E} + \frac{R^2}{(\epsilon R)^{1+\alpha}}\|f\|_{0,F}, \tag{4.15}$$

$$\|\hat{P}f\|_{2,E} \lesssim \|f\|_{0,E} + \frac{R^2}{(\epsilon R)^{1+\alpha}}\|f\|_{0,F}, \text{ and}$$

$$\|\hat{P}f\|_{2,F} \lesssim \frac{1}{(\epsilon R)^{\alpha}} \frac{1}{\epsilon^2 R^{2+\delta}}\left(\|f\|_{0,E} + \frac{R^2}{(\epsilon R)^{1+\alpha}}\|f\|_{0,F}\right)$$

## 5 - Existence and Embeddedness.

**5.1 - The Schauder fixed-point theorem.**   At this point, it is convenient to modify slightly the norms introduced in (4.7). We define,

$$\|f\|'_{m,F,\text{Höl}} := \|f|_{A(2R,\infty)}\|_{C^{m,\alpha}_{\gamma,\epsilon}(\mathbb{R}^2)},$$
$$\|f\|'_{m,F,\text{Sob}} := \|f|_{A(2R,\infty)}\|_{H^{m}_{\gamma,\epsilon}(\mathbb{R}^2)}, \text{ and} \tag{5.1}$$
$$\|f\|'_{m,F} := \|f\|'_{m,F,\text{Höl}} + \frac{1}{(\epsilon R)}\|f\|'_{m,F,\text{Sob}}.$$

Observe that, by (4.12), this does not affect (4.15). In addition, we will also ignore the factor of $1/\psi^\epsilon$ used in the definitions of the operators $\hat{J}^\epsilon_{\hat{\Sigma}}$, $X^\epsilon_{\hat{\Sigma}}$, $Y^\epsilon_{\hat{\Sigma}}$ and $Z^\epsilon_{\hat{\Sigma}}$. Indeed, it is readily shown that the operator of multiplication by this function is uniformly bounded, independent of $\Lambda$, with respect to the norms $\|\cdot\|_{0,E}$ and $\|\cdot\|_{0,F}$ and therefore also does not affect these estimates.

For all non-negative, integer $m$, for all $\alpha \in [0,1]$ and for all real $\gamma$, let $\mathcal{L}^{m,\alpha}_\gamma(\hat{\Sigma})$ denote the Frechet space of $m$-times differentiable functions $f : \hat{\Sigma} \to \mathbb{R}$ which are invariant under all horizontal symmetries of the Costa-Hoffman-Meeks surface and which satisfy

$$\|f\|_{m,E}, \|f\|'_{m,F} < \infty.$$

Now let $\tilde{\mathcal{E}} : \mathbb{R}^3 \oplus \mathbb{R}^3 \oplus \mathbb{R}^3 \oplus C^\infty_0(\hat{\Sigma}) \to C^\infty(\hat{\Sigma}, \mathbb{R}^3)$ be the family of immersions constructed in Section 3.2 and let $\mathcal{H} : \mathbb{R}^3 \oplus \mathbb{R}^3 \oplus \mathbb{R}^3 \oplus \mathcal{L}^{2,\alpha}_\gamma(\hat{\Sigma}) \to \mathcal{L}^{0,\alpha}_\gamma(\hat{\Sigma})$ be the corresponding mean curvature functional with respect to the metric $\overline{g}^\epsilon$.

**Lemma 5.1.1**

$$\|\mathcal{H}(0,0,0,0)\|_{0,E} \lesssim R^{\delta-2}, \text{ and}$$
$$\|\mathcal{H}(0,0,0,0)\|'_{0,F} = 0. \tag{5.2}$$

**Proof:** Denote $\psi := \mathcal{H}(0,0,0,0)$. We study this function over the regions $\Sigma \cap (B(R_0) \times \mathbb{R})$, $\Sigma \cap (A(R_0, R) \times \mathbb{R})$, $\Sigma \cap (A(R, 2R) \times \mathbb{R})$ and $\Sigma \cap (A(2R, \infty) \times \mathbb{R})$. First, over the cylinder $B(R_0) \times \mathbb{R}$,

$$\left\|(\overline{g}^\epsilon - \overline{g}^0)|_{B(R_0) \times \mathbb{R}}\right\|_{C^{1,\alpha}} \lesssim \epsilon^2.$$

Thus, since $\mathcal{H}$ is a smooth functional of the metric and the surface, and since $\hat{\Sigma} \cap (B(R_0) \times \mathbb{R})$ is minimal with respect to $\overline{g}^0$, it follows that

$$\left\|\psi|_{\hat{\Sigma} \cap (B(R_0) \times \mathbb{R})}\right\|_{C^{0,\alpha}} \lesssim \epsilon^2 \lesssim R^{\delta-2}. \tag{5.3}$$

Consider now one of the components of $\hat{\Sigma} \cap (A(R_0, \infty) \times \mathbb{R})$ and let $u : A(R_0, \infty) \to \mathbb{R}$ be its profile. Observe that the mean curvature of the graph of $u$ with respect to the metric $\overline{g}$ takes the form

$$\mathcal{H}(u, \overline{g}) = F(\overline{g}, Du)^{ij} \text{Hess}(u)_{ij} + G(\overline{g}, D\overline{g}, Du)^i u_i,$$

where $F$ and $G$ are smooth functions of their arguments and $G(\overline{g}, 0, Du)^i = 0$. However, over the region $A(R_0, R) \times \mathbb{R}$,

$$\overline{g}^\epsilon - \overline{g}^0 = \mathrm{O}\left(\frac{\epsilon^2 r^2}{r^k}\right), \text{ and}$$

$$u = a + c\mathrm{Log}(r) + \mathrm{O}\left(\frac{1}{r^{2+k}}\right),$$

Since $\hat{\Sigma} \cap (A(R_0, R) \times \mathbb{R})$ is minimal with respect to $\overline{g}^0$, it follows that, over this region,

$$\psi = \mathrm{O}\left(\frac{\epsilon^2}{r^k}\right),$$

so that

$$\left\| \psi|_{\hat{\Sigma} \cap (A(R_0, R) \times \mathbb{R})} \right\|_{C^{0,\alpha}_{\delta+2, \mathrm{SF}}} \lesssim \epsilon^2 R^{\delta+2} \lesssim R^{\delta-2}. \tag{5.4}$$

Now let $u^\epsilon$ be the profile of the $\epsilon$-hyperbolic end with constant term $a$ and logarithmic parameter $c$. Over $A(R, 2R)$,

$$(u - u^\epsilon) = \mathrm{O}\left(\frac{1}{r^{2+k}}\right).$$

Thus, if $\hat{J}^\epsilon$ denotes the modified Jacobi operator of $u^\epsilon$ then, by (2.1) and (3.11), over this annulus,

$$\hat{J}^\epsilon(u - u^\epsilon) = \mathrm{O}\left(\frac{1}{r^{4+k}}\right).$$

Since the mean curvature is a smooth functional of the profile, and since $u^\epsilon$ is minimal with respect to $\overline{g}^\epsilon$, it follows that $\psi$ differs from $\hat{J}^\epsilon(u - u^\epsilon)$ by a term which is quadratic in $(u - u^\epsilon)$. Thus, over $\hat{\Sigma} \cap (A(R, 2R) \times \mathbb{R})$,

$$\psi = \mathrm{O}\left(\frac{1}{r^{4+k}}\right),$$

so that

$$\left\| \psi|_{\hat{\Sigma} \cap (A(R, 2R) \times \mathbb{R})} \right\|_{C^{0,\alpha}_{\delta+2, \mathrm{SF}}} \lesssim R^{\delta-2}. \tag{5.5}$$

Finally, by construction, $\psi$ vanishes over $\hat{\Sigma} \cap (A(2R, \infty) \times \mathbb{R})$ and the result now follows upon combining (5.3), (5.4) and (5.5). $\square$

For all $(u, v)$, let $\hat{J}^{u,v}$, $X^{u,v}$, $Y^{u,v}$ and $Z^{u,v}$ denote respectively the modified Jacobi operator and the perturbation operators of the surface $\Sigma^{u,v} := \hat{\mathcal{E}}(u, v, 0, 0)(\Sigma)$ defined in Section 4.2, above. Likewise, let $(\hat{P}^{u,v}, \hat{U}^{u,v}, \hat{V}^{u,v}, \hat{W}^{u,v})$ be the right inverse of $(\hat{J}^{u,v}, X^{u,v}, Y^{u,v}, Z^{u,v})$. Observe now that the estimates derived in the previous sections also apply to these operators uniformly for $(u, v)$ in a small neighbourhood of the origin.

**Lemma 5.1.2**

*If $\|u\|$, $\|v\|$, $\|w\|$ and $\|f\|_{2,F}$ are sufficiently small, then*

$$\|\mathcal{H}(u,v,w,f) - \mathcal{H}(0,0,0,0) - \hat{J}^{u,v}f - X^{u,v}u - Y^{u,v}v\|_{0,E} \lesssim \|u\|^2 + \|v\|^2 + \|f\|_{2,E}^2. \quad (5.6)$$

**Proof:** First observe that, over $\hat{\Sigma} \cap (B(4R) \times \mathbb{R})$,

$$\mathcal{H}(u,v,w,f) = \mathcal{H}(u,v,0,f).$$

Next, by Taylor's Theorem

$$\|\mathcal{H}(u,v,0,f) - \mathcal{H}(u,v,0,0) - \hat{J}^{u,v}f\|_{0,E} \lesssim \|f\|_{2,E}^2, \text{ and}$$

$$\|\mathcal{H}(u,v,0,0) - \mathcal{H}(0,0,0,0) - X^{u,v}u - Y^{u,v}v\|_{0,E} \lesssim \|u\|^2 + \|v\|^2.$$

The result follows upon combining these relations. $\square$

**Lemma 5.1.3**

*For sufficiently small $\alpha \in [0,1]$, for sufficiently small $\gamma$ and for sufficiently large $\Lambda$, if*

$$\|f\|_{2,G}' < \frac{1}{(\epsilon R)^{2-\alpha}},$$

*then*

$$\|\mathcal{H}(u,v,w,f) - \mathcal{H}(0,0,0,0) - \hat{J}^{u,v}f - Z^{u,v}w\|_{0,E}'$$
$$\lesssim \frac{\epsilon^2}{R(\epsilon R)^{2\alpha}}\|w\|^2 + \epsilon^3(\epsilon R)^{1-2\alpha}\left(\|f\|_{2,F}'\right)^2. \quad (5.7)$$

**Proof:** First observe that, over $\hat{\Sigma} \cap (A(2R,\infty) \times \mathbb{R})$,

$$\mathcal{H}(u,v,0,0) = \mathcal{H}(0,0,0,0) = 0.$$

Next, rescaling (2.31) yields

$$\left\|\mathcal{H}(u,v,0,f) - \mathcal{H}(u,v,0,0) - \hat{J}^{u,v}_{\hat{\Sigma}}f\right\|_{0,F}' \lesssim \epsilon^3(\epsilon R)^{1-2\alpha}\left(\|f\|_{2,F}'\right)^2. \quad (5.8)$$

Finally, observe that, over $\hat{\Sigma} \setminus (A(1/2\epsilon, 1/\epsilon) \times \mathbb{R})$,

$$\mathcal{H}(u,v,w,f) = \mathcal{H}(u,v,0,f),$$

whilst, over $\hat{\Sigma} \cap (A(1/2\epsilon, 1/\epsilon) \times \mathbb{R})$,

$$\mathcal{H}(u,v,w,f) = \mathcal{H}\left(u,v,0,f + \sum_{i=1}^{3} b_i \mathbb{I}_i(1 - \chi_\epsilon')\right), \quad (5.9)$$

where, for $1 \leq i \leq 3$, $\mathbb{I}_i$ is the indicator function of the $i$'th component of $\hat{\Sigma} \cap (A(1/2\epsilon, 1/\epsilon) \times \mathbb{R})$. However,

$$Z^\epsilon_{\hat{\Sigma}} w = -\sum_{i=1}^{3} b_i \mathbb{I}_i \hat{J}^\epsilon_{\hat{\Sigma}} \chi_\epsilon'.$$

Furthermore, for each $i$,

$$\left\|\mathbb{I}_i(1 - \chi_\epsilon')|_{\hat{\Sigma} \cap (A(1/2\epsilon, 1/\epsilon) \times \mathbb{R})}\right\|_{2,F} \lesssim \frac{1}{(\epsilon R)}. \quad (5.10)$$

The result now follows upon combining (5.8), (5.9) and (5.10). $\square$

**Theorem 5.1.4**

*For all $\delta \in ]1, 2[$, for all sufficiently small $\alpha$, for all sufficiently small $\gamma$, and for all sufficiently large $\Lambda$, there exist $u$, $v$, $w$ and $f$ such that*

$$\mathcal{H}(u, v, w, f) = 0.$$

*Furthermore,*

$$\|u\|, \|v\|, \|w\|, \|f\|_{2,E} \lesssim R^{\delta-2}, \ \|f\|_{2,F} \lesssim \frac{1}{(\epsilon R)^\alpha \epsilon^2 R^4}. \tag{5.11}$$

**Proof:** Fix $\delta \in ]1, 2[$ and $\alpha$ and $\gamma$ small. Set $\psi_0 := \mathcal{H}(0, 0, 0, 0)$. By (5.2),

$$\|\psi_0\|_{0,E} \lesssim R^{\delta-2}, \ \|\psi_0\|'_{0,F} = 0.$$

It follows by (2.1) and (4.15) that for all $(u, v)$ sufficiently close to $(0, 0)$,

$$\|\hat{U}^{u,v}\psi_0\|, \|\hat{V}^{u,v}\psi_0\|, \|\hat{W}^{u,v}\psi_0\|, \|\hat{P}^{u,v}\psi_0\|_{2,E} \leq BR^{\delta-2}, \ \|\hat{P}^{u,v}\psi_0\|'_{2,F} \leq \frac{B}{(\epsilon R)^\alpha \epsilon^2 R^4},$$

for some constant $B > 0$, independent of $\Lambda$. Let $\Omega \subseteq \mathbb{R}^3 \oplus \mathbb{R}^3 \oplus \mathbb{R}^3 \oplus \mathcal{L}^{2,\alpha}_\gamma(\hat{\Sigma})$ be the set of all quadruplets $(u, v, w, f)$ such that

$$\|u\|, \|v\|, \|w\|, \|f\|_{2,E} \leq 2BR^{\delta-2}, \ \|f\|'_{2,F} \leq \frac{2B}{(\epsilon R)^\alpha \epsilon^2 R^4}.$$

Observe that $\Omega$ is convex and, by the Arzela-Ascoli Theorem, for all $\alpha' < \alpha$ and $\gamma' < \gamma$, $\Omega$ is a compact subset of $\mathbb{R}^3 \oplus \mathbb{R}^3 \oplus \mathbb{R}^3 \oplus \mathcal{L}^{2,\alpha'}_{\gamma'}(\hat{\Sigma})$. For $\phi := (u, v, w, f) \in \Omega$, define

$$\Phi(\phi) := -(\hat{U}^{u,v}(\psi_0 + \psi), \hat{V}^{u,v}(\psi_0 + \psi), \hat{W}^{u,v}(\psi_0 + \psi), \hat{P}^{u,v}(\psi_0 + \psi)),$$

where

$$\psi := \mathcal{H}(u, v, w, f) - \mathcal{H}(0, 0, 0, 0) - \hat{J}^{u,v}f - X^{u,v}u - Y^{u,v}v - Z^{u,v}w.$$

By (2.1), (5.6) and (5.7),

$$\|\psi\|_{0,E} \lesssim R^{2\delta-4}, \ \|\psi\|'_{0,F} \lesssim \frac{1}{(\epsilon R)^{4\alpha} R^7},$$

so that, by (2.1) and (4.15), for sufficiently large $\Lambda$, $\Phi$ maps $\Omega$ to itself. Furthermore, for all $\alpha' < \alpha$ and $\gamma' < \gamma$, $\Phi$ is continuous with respect to the topology of $\mathcal{L}^{2,\alpha'}_{\gamma'}(\hat{\Sigma})$. It follows by the Schauder fixed point theorem (c.f. [5]) that $\Phi$ has a fixed point $\phi = (u, v, w, f)$ in $\Omega$. However, for such a $\phi$,

$$\begin{aligned}
\mathcal{H}(u, v, w, f) &= \mathcal{H}(0, 0, 0, 0) + \psi + \hat{J}^{u,v}f + X^{u,v}u + Y^{u,v}v + Z^{u,v}w \\
&= (\psi_0 + \psi) - \hat{J}^{u,v}\hat{P}^{u,v}(\psi_0 + \psi) - X^{u,v}\hat{U}^{u,v}(\psi_0 + \psi) \\
&\quad - Y^{u,v}\hat{V}^{u,v}(\psi_0 + \psi) - Z^{u,v}\hat{W}^{u,v}(\psi_0 + \psi) \\
&= (\psi_0 + \psi) - (\psi_0 + \psi) \\
&= 0,
\end{aligned}$$

as desired. $\square$

**Theorem 5.1.5**

*Let $(u, v, w, f)$ be as in Theorem 5.1.4. For sufficiently large $\Lambda$, the surface $\tilde{\mathcal{E}}(u, v, w, f)$ is embedded.*

**Proof:** Recall that the joined surface is denoted by $\hat{\Sigma}$. Let $S$ denote the image of $\tilde{\mathcal{E}}(u, v, w, f)$. We rescale both $\hat{\Sigma}$ and $S$ to obtain immersed surfaces in $\mathbb{H}^3$. Observe now that $\hat{\Sigma} \cap (A(2\epsilon R, \infty) \times \mathbb{R})$ consists of 3 distinct hyperbolic minimal ends which we denote by $E_+$, $E_0$ and $E_-$ respectively. Let $u_+$, $u_0$ and $u_-$ be the respective profiles of these ends, and let $v_+$, $v_0$ and $v_-$ be their derivatives in the radial direction. Observe that

$$u_+(\epsilon R) > u_0(\epsilon R) > u_-(\epsilon R),$$

and

$$v_+(\epsilon R) > v_0(\epsilon R) > v_-(\epsilon R).$$

Since $v_+$, $v_0$ and $v_-$ are all solutions of the same first order ODE, it follows that $v_+(r) > v_0(r) > v_-(r)$ for all $r$. In particular, the ends $E_+$, $E_0$ and $E_+$ are separated vertically by a (euclidean) distance of no less than $\eta$, where $\eta \sim \epsilon \text{Log}(R)$. Let $\Omega_+$, $\Omega_0$ and $\Omega_-$ denote the open sets of points lying at a vertical (euclidean) distance of no more than $\eta/2$ from $E_+$, $E_0$ and $E_-$ respectively. In particular, these 3 sets are disjoint.

Now let $E'_+$, $E'_0$ and $E'_-$ be the three ends of $S$. Over the annulus $A(\epsilon R, 2\epsilon R)$, by (5.11),

$$\left\| \epsilon f|_{A(R, 2R)} \right\|_{C^0} \lesssim \epsilon R^{-\delta} \left\| f|_{(A(R, 2R))} \right\|_{2, F} \lesssim \epsilon R^{-2},$$

so that, over this annulus, $E'_+$ lies strictly above $E'_0$, and $E'_0$ lies strictly above $E'_-$. However, by (2.29) and (5.11) again,

$$\left\| \epsilon f \right\|'_{1, F, \text{Höl}} \lesssim \frac{1}{(\epsilon R)^{3\alpha} R^3}.$$

Bearing in mind the definition of the norm $\| \cdot \|_{1, G, H}$, it follows that for sufficiently large $\Lambda$, $E'_+$, $E'_0$ and $E'_-$ are all graphs over $A(\epsilon R, \infty)$. Furthermore, for some large $R'$, the intersections of each of these ends with $A(R', \infty) \times \mathbb{R}$ are contained in $\Omega_+$, $\Omega_0$ and $\Omega_-$ respectively. In particular, outside $B(R') \times \mathbb{R}$, $E'_+$ lies strictly above $E'_0$ and $E'_0$ lies strictly above $E'_-$. Since minimality is preserved by vertical translations in $\mathbb{H}^3$, it follows by the strong maximum principle that, over the whole of $A(\epsilon R, \infty)$, $E'_+$ lies strictly above $E'_0$ and $E'_0$ lies strictly above $E_-$. This completes the proof. $\square$

# A - Terminology, conventions and standard results.

## A.1 - General definitions.

(1) $\mathbb{R}^2$ and $\mathbb{R}^3$ will denote respectively two- and three-dimensional euclidian spaces.

(2) $\mathbb{R}^2$ will be identified with the $x - y$ plane in $\mathbb{R}^3$.

(3) $e_x$, $e_y$ and $e_z$ will denote the vectors of a canonical basis of $\mathbb{R}^3$.

(4) $D$ will denote the total differentiation operator over $\mathbb{R}^2$.

(5) $r$ will denote the radial distance to the origin in $\mathbb{R}^2$ as well as the radial distance to the $z$-axis in $\mathbb{R}^3$.

(6) For all $R$, $C(R)$ will denote the circle of radius $R$ about the origin in $\mathbb{R}^2$.

(7) For all $R$, $B(R)$ will denote the closed disk of radius $R$ about the origin in $\mathbb{R}^2$

(8) For all $R_1 < R_2$, $A(R_1, R_2)$ will denote the closed annulus of inner radius $R_1$ and outer radius $R_2$ about the origin in $\mathbb{R}^2$.

(9) Let $\chi : [0, \infty[ \to \mathbb{R}$ be a non-negative, non-increasing function such that $\chi = 1$ over $[0, 1]$ and $\chi = 0$ over $[2, \infty[$. For all $R > 0$, define $\chi_R : \mathbb{R}^2 \to \mathbb{R}$ by

$$\chi_R(x) := \chi\left(\frac{r}{R}\right). \tag{A.1}$$

We call $\chi_R$ the *cut-off function* of the *transition region* $A(R, 2R)$.

## A.2 - A-priori estimates.

(1) Given two variable quantities $a$ and $b$, we will write

$$a \lesssim b. \tag{A.2}$$

to mean that there exists a constant $C$, which for the purposes of this paper will be considered universal, such that

$$a \leq Cb.$$

(2) Given a function $f$ and a sequence of functions $(g_m)$, we will write

$$f = \mathrm{O}(g_m). \tag{A.3}$$

to mean that there exists a sequence $(C_m)$ of constants, which for the purposes of this paper will be considered universal, such that the relation

$$|D^m f| \leq C_m g_m$$

holds pointwise for all $m$. The indexing variable of the sequence $(g_m)$ should be clear from the context. In certain cases, every element of the sequence $(g_m)$ may be the same. It should also be clear from the context when this is the case.

## A.3 - General surface geometry.
Let $\Omega$ be a domain in $\mathbb{R}^2$. Let $u : \Omega \to \mathbb{R}$ be a smooth function. Let $\Sigma^u$ be its graph in $\mathbb{R}^3$. Let $\overline{g}$ be a smooth metric over $\mathbb{R}^3$. The following geometric objects will be constructed with respect to this metric. Furthermore, all objects defined over $\Sigma^u$ will be identified with objects defined over $\Omega$ via the canonical projection.

(1) $N^u$ will denote the unit normal vector field over $\Sigma^u$.

(2) $\pi^u$ will denote the orthogonal projection onto the tangent space of $\Sigma^u$.

(3) $g^u$ will denote the intrinsic metric of $\Sigma^u$.

(4) $\nabla^u$ will denote the Levi-Civita covariant derivative of $\Sigma^u$.

(5) $\mathrm{Hess}^u$ will denote the Hessian operator of $\Sigma^u$.

(6) $\Delta^u$ will denote the Laplace-Beltrami operator of $\Sigma^u$.

(7) $\mathrm{II}^u$ will denote the second fundamental form of $\Sigma^u$ with respect to the normal direction $N^u$.

(8) $A^u$ will denote the shape operator of $\Sigma^u$ with respect to the normal direction $N^u$.

(9) $H^u$ will denote the mean curvature - that is, the *sum* of the principal curvatures - of $\Sigma^u$ with respect to the normal direction $N^u$

(10) $J^u$ will denote the Jacobi operator of $\Sigma^u$. Recall that $J^u$ is defined as follows. First the function $\mathcal{E} : C_0^\infty(\Sigma^u) \times \Sigma^u \to \mathbb{R}^3$ is defined by

$$\mathcal{E}(f, x) := x + f(x) N^u(x).$$

For all sufficiently small $f$, $\mathcal{E}(f, \cdot)$ is an embedding. The function $\mathcal{H} : C_0^\infty(\Sigma^u) \times \Sigma^u \to \mathbb{R}$ is then defined such that, for all suitable $f$ and for all $x$, $\mathcal{H}(f, x)$ is the mean curvature of $\mathcal{E}(f, \cdot)$ at the point $x$. The Jacobi operator is now defined by

$$(J^u f)(x) := \left. \frac{\partial}{\partial t} \mathcal{H}(tf, x) \right|_{t=0}. \tag{A.4}$$

Recall that $J^u$ is given explicitly by

$$J^u f = -\big[\mathrm{Ric}\big(N^u, N^u\big) + \mathrm{Tr}\big(A^u A^u\big)\big] f - \Delta^u f,$$

where Ric here denotes the Ricci curvature tensor of $\bar{g}$.

The following elementary relations will also prove useful.

(11) For any point $x \in \Sigma^u$, and for any function $f$ defined over a neighbourhood of $x$ in $\mathbb{R}^3$,

$$\nabla^u(f)(x) = Df(x) - \langle Df(x), N^u(x)\rangle N^u(x). \tag{A.5}$$

(12) For any point $x \in \Sigma^u$, for any function $f$ defined over a neighbourhood of $x$ in $\mathbb{R}^3$, and for all vectors $X$ and $Y$ tangent to $\Sigma^u$ at $x$,

$$\mathrm{Hess}^u(f)(x)(X, Y) = D^2 f(x)(X, Y) - \langle Df(x), N^u(x)\rangle \mathrm{II}^u(x)(X, Y). \tag{A.6}$$

(13) If $J^u \phi = 0$ and if $M_\phi$ denotes the operator of multiplication by $\phi$, then

$$(M_\phi^{-1} J^u M_\phi) f = -\Delta^u f - \frac{2}{\phi} g^{u,ij} \phi_i f_j. \tag{A.7}$$

**A.4 - Minimal graphs in hyperbolic space.** We continue to use to the notation of Section A.3. Let $\bar{g}$ and $g$ be the hyperbolic metrics defined in Section 2.1.

(1) We define

$$\mu^u := \frac{\cosh(r)}{\sqrt{1 + \cosh^2(r)\|\nabla^g u\|_g^2}}. \tag{A.8}$$

(2) The upward pointing unit normal vector field over $\Sigma^u$ is given by

$$N^u = \mu^u \left( -\nabla^g u, \frac{1}{\cosh^2(r)} \right), \tag{A.9}$$

where $\nabla^g$ here denotes the Levi-Civita covariant derivative of the metric $g$.

(3) The intrinsic metric of $\Sigma^u$ is given by

$$g_{ij}^u = g_{ij} + \cosh^2(r)u_i u_j. \tag{A.10}$$

(4) The inverse of this metric is given by

$$g^{u,ij} = g^{ij} - (\mu^u)^2 u^i u^j, \tag{A.11}$$

where the indices are raised with respect to $g$.

(5) The relative Christoffel symbol of $g^u$ with respect to $g$ is defined by

$$\Gamma_{ij}^k := \left( \nabla_{e_i}^u e_j - \nabla_{e_i}^g e_j \right)^k. \tag{A.12}$$

By Koszul's formula, it is given by

$$\Gamma_{ij}^k = g^{u,kp}\cosh^2(r)\text{Hess}^g(u)_{ij}u_p + g^{u,kp}\cosh(r)\sinh(r)\left( u_i u_p r_j + u_j u_p r_i - u_i u_j r_p \right), \tag{A.13}$$

where $\text{Hess}^g$ here denotes the Hessian operator of the metric $g$.

(6) The Hessian operator of $\Sigma^u$ is given by

$$\text{Hess}^u(f)_{ij} := \text{Hess}^g(f)_{ij} - \Gamma_{ij}^k f_k. \tag{A.14}$$

(7) The Laplace-Beltrami operator of $\Sigma^u$ is given by

$$\begin{aligned}
\Delta^u f = &\ g^{u,ij}\text{Hess}^g(f)_{ij} - \cosh^2(r)g^{u,ij}g^{u,kp}\text{Hess}^g(u)_{ij}u_p f_k \\
&- 2\cosh(r)\sinh(r)g^{u,ij}g^{u,kp}u_i u_p r_j f_k \\
&+ \cosh(r)\sinh(r)g^{u,ij}g^{u,kp}u_i u_j r_p f_k.
\end{aligned} \tag{A.15}$$

(8) The second fundamental form of $\Sigma^u$ is given by

$$\begin{aligned}
\Pi_{ij}^u = &-\mu^u\text{Hess}^g(u)_{ij} - \mu^u\tanh(r)\left( r_i u_j + u_i r_j \right) \\
&- \mu^u\sinh(r)\cosh(r)u_r u_i u_j.
\end{aligned} \tag{A.16}$$

(9) The mean curvature of $\Sigma^u$ is given by

$$\begin{aligned}
H^u = &-\mu^u g^{u,ij}\text{Hess}^g(u)_{ij} - 2\mu^u\tanh(r)g^{u,ij}r_i u_j \\
&- \mu^u\sinh(r)\cosh(r)g^{u,ij}u_r u_i u_j.
\end{aligned} \tag{A.17}$$

**A.5 - Function spaces.** Let $X$ be a complete riemannian manifold.

(1) For all $\alpha \in [0,1]$, the *Hölder semi-norm* of order $\alpha$ is defined over functions over $X$ by

$$[f]_\alpha := \underset{x \neq y \in X}{\text{Sup}} \frac{|f(x) - f(y)|}{d(x,y)^\alpha}. \tag{A.18}$$

In particular $[f]_0$ is the *total variation* of $f$ and $[f]_1$ is its Lipschitz semi-norm.

(2) For all $\alpha \in [0,1]$, we define

$$\delta^\alpha f(x) := [f|_{B_1(x)}]_\alpha. \tag{A.19}$$

(3) For all non-negative integer $m$ and for all $\alpha \in [0,1]$, the Hölder norm of order $(m, \alpha)$ is defined over $k$-times differentiable functions over $X$ by

$$\|f\|_{C^{m,\alpha}(X)} := \sum_{k=0}^{m} \|D^k f\|_{C^0(X)} + \|\delta^\alpha D^m f\|_{C^0([0,\infty[)}. \tag{A.20}$$

The Hölder space of order $(m, \alpha)$ over $X$ will be denoted by $C^{m,\alpha}(X)$. This is the space of all $m$-times differentiable functions $f : X \to \mathbb{R}$ such that $\|f\|_{C^{k,\alpha}(X)} < \infty$.

(4) For all $p \in [1, \infty[$, the Lebesgue norm with exponent $p$ is defined over measurable functions over $X$ by

$$\|f\|_{L^p(X)}^p := \int_X |f|^p \, \mathrm{dVol}. \tag{A.21}$$

The Lebesgue space with exponent $p$ will be denoted by $L^p(X)$. This is the space of all measurable functions $f : X \to \mathbb{R}$ such that $\|f\|_{L^p} < \infty$.

(5) For all non-negative integer $m$, the Sobolev norm of order $m$ is defined over $m$-times locally $L^2$ differentiable distributions over $X$ by

$$\|f\|_{H^m(X)} := \sum_{k=0}^{m} \|D^i f\|_{L^2(X)}. \tag{A.22}$$

The Sobolev space of order $m$ over $X$ will be denoted by $H^m(X)$. This is the space of distributions $f$ over $X$ all of whose derivatives up to and including order $m$ are elements of $L^2(X)$.

We will also make use of the following relations.

(6) For all $\alpha \in [0,1]$, and for all $f : X \to \mathbb{R}$,

$$[f]_\alpha \leq [f]_0^{1-\alpha}[f]_1^\alpha \leq 2^{1-\alpha}\|f\|_{C^0(X)}^{1-\alpha}[f]_1^\alpha. \tag{A.23}$$

(7) For all $\alpha \in [0,1[$, for all $\beta \in ]0,1]$, and for all differentiable $f : X \to \mathbb{R}$,

$$\|Df\|_{C^0} \leq 2[f]_\alpha^{\frac{\beta}{1+(\beta-\alpha)}}[Df]_\beta^{\frac{1-\alpha}{1+(\beta-\alpha)}}. \tag{A.24}$$

(8) For all $\alpha \in [0,1]$ and for all $f, g : X \to \mathbb{R}$,

$$[fg]_\alpha \leq \|f\|_{C^0}[g]_\alpha + [f]_\alpha \|g\|_{C^0(X)}. \tag{A.25}$$

(9) If $\Phi : \mathbb{R}^n \to \mathbb{R}^m$ is a smooth function and $\|u\|_{C^0} \leq B$, then there exists a constant $C$ which only depends on $\Phi$ and $B$ such that

$$\begin{aligned} \|\Phi(u)\|_{C^0} &\leq C, \text{ and} \\ [\Phi(u)]_\alpha &\leq C[u]_\alpha. \end{aligned} \tag{A.26}$$

(10) If $\Phi : \mathbb{R}^n \to \mathbb{R}^m$ is a smooth function and $\|u\|_{C^0}, \|v\|_{C^0} \leq B$, then there exists a constant $C$ which only depends on $\Phi$ and $B$ such that

$$\|\Phi(u) - \Phi(v)\|_{C^{0,\alpha}} \leq C(1 + [u]_\alpha + [v]_\alpha)\|u - v\|_{C^{0,\alpha}}. \tag{A.27}$$

(11) We readily verify that all the surfaces studied in this paper are sufficiently regular at infinity for the Sobolev embedding theorem to hold. That is, for all non-negative integer $n$, for all $m + \alpha < n - 1$, and for all suitable $f : X \to \mathbb{R}$,

$$\|f\|_{C^{m,\alpha}} \lesssim \|f\|_{H^n}. \tag{A.28}$$

# B - Bibliography.